\documentclass[10pt]{amsart}
\allowdisplaybreaks[4]
\allowdisplaybreaks[4]
\usepackage{amsthm}
\usepackage{graphicx}

\usepackage{amssymb,color}

\definecolor{c20}{rgb}{0.,0.7,0.}
\definecolor{c30}{rgb}{0.,0.,1.}
\definecolor{c40}{rgb}{1,0.1,0.7}
\definecolor{c50}{rgb}{1,0,0}
\definecolor{c60}{rgb}{0,0.9,0.1}

\usepackage{amssymb,color}

\newcommand{\kb}[1]{\boldsymbol{#1}}
\newcommand{\vk}[1]{\kb{#1}}

\newcommand{\E}[1]{\mathbb{E}\left(#1\right)}

\newcommand{\pl}[1]{\mathbb{P}\left(#1 \right)}

\newcommand{\limit}[1]{\lim_{#1 \to   \infty}}

\newcommand{\BQN}{\begin{eqnarray}}
\newcommand{\EQN}{\end{eqnarray}}
\newcommand{\BQNY}{\begin{eqnarray*}}
\newcommand{\EQNY}{\end{eqnarray*}}

\newcommand{\BS}{\begin{sat}}
\newcommand{\ES}{\end{sat}}
\newcommand{\BT}{\begin{theo}}
\newcommand{\ET}{\end{theo}}
\newcommand{\BK}{\begin{korr}}
\newcommand{\EK}{\end{korr}}

\newcommand{\BD}{\begin{de}}
\newcommand{\ED}{\end{de}}

\newcommand{\BIT}{\begin{itemize}}
\newcommand{\EIT}{\end{itemize}}
\newcommand{\BDI}{\begin{description}}
\newcommand{\EDI}{\end{description}}

\newcommand{\BRM}{\begin{remarks}}
\newcommand{\ERM}{\end{remarks}}

\newcommand{\BEL}{\begin{lem}}
\newcommand{\EEL}{\end{lem}}

\newtheorem{theo}{Theorem}[section]
\newtheorem{sat}[theo]{Proposition}
\newtheorem{de}[theo]{Definition}
\newtheorem{lem}[theo]{Lemma}

\newtheorem{korr}[theo]{Corollary}
\newtheorem{remark}[theo]{Remark}
\newtheorem{remarks}[theo]{Remarks}

\newcommand{\prooflem}[1]{\textbf{Proof of Lemma} \ref{#1}}

\newcommand{\COM}[1]{}

\newcommand{\QED}{\hfill $\Box$}

\def\rw{\rightarrow}
\def\RW{\Rightarrow}

\def\IF{\infty}

\date{}

\begin{document}

\title{  Extremes of stationary Gaussian storage models}

\author{Krzysztof D\c{e}bicki}
\address{Krzysztof D\c{e}bicki, Mathematical Institute, University of Wroc\l aw, pl. Grunwaldzki 2/4, 50-384 Wroc\l aw, Poland}
\email{Krzysztof.Debicki@math.uni.wroc.pl}

\author{Peng Liu}
\address{School of Mathematical Sciences and LPMC, Nankai University, Tianjin 300071, China and Department of Actuarial Science, University of Lausanne, UNIL-Dorigny 1015 Lausanne, Switzerland}
\email{peng.liu@unil.ch}

\bigskip

\date{\today}
 \maketitle

\bigskip
{\bf Abstract:}
For the stationary storage process $\{Q(t), t\ge0\}$, with
$
Q(t)=\sup_{ s \ge t}\left(X(s)-X(t)-c(s-t)^\beta\right),
$
where $\{X(t),t\ge 0\}$ is a centered Gaussian process with stationary increments, $c>0$
and $\beta>0$ is chosen such that $Q(t)$ is finite a.s.,
we derive exact asymptotics of $\pl{\sup_{t\in [0,T_u]} Q(t)>u}$
and
$\pl{\inf_{t\in [0,T_u]} Q(t)>u}$, as $u\to\infty$.
As a by-product we find conditions under which {\it strong Piterbarg property} holds.




\section{Introduction}

Let $\{X(t),t\ge 0\}$ be a centered Gaussian process with stationary increments,
a.s. continuous sample paths and variance function $\sigma^2(t)$.
Given $c>0$ and $\beta>0$ , consider the stationary storage process $\{Q(t), t\ge0\}$, with
\BQN\label{sto}
Q(t)=\sup_{ s \ge t}\left(X(s)-X(t)-c(s-t)^\beta\right),\ \ \ t\geq 0,
\EQN
where $c>0$ and $\beta>0$ is chosen appropriately to guarantee a.s. finiteness of $Q(t)$.

The stimulus to analyze distributional properties of $\{Q(t), t\ge0\}$
stems, for instance,  from its straightforward relation with the theory of {\it reflected}
Gaussian processes, its applications in widely investigated
Gaussian fluid queueing models and,
by duality, its importance in risk theory.
In particular, for $\beta=1$, by Reich representation \cite{REICH},
$Q(t)$ describes the stationary amount of substance in reservoir, where
the inflow to the reservoir in time interval $[s,t]$ equals to
$X(t)-X(s)$ and the rate of outflow is $c$.

Motivated by the above applications,
$Q(0)$ has been studied in the literature under different levels of generality,
e.g., \cite{Norros94}, \cite{HP99}, \cite{DE2002}, \cite{HP2004},
\cite{DI2005}, \cite{HA2013}, \cite{PengLith}.
Particularly vast interest has been paid to the analysis of storage models, where
$X(t)=B_H(t)$ is a fractional Brownian motion (fBm) with Hurst index $H\in (0,1)$
and $\beta=1$, leading
to derivation of exact asymptotics of $\pl{Q(0)>u}$ as $u\to\infty$ in \cite{HP99}
and a surprising asymptotic equivalence
\begin{eqnarray}\label{strong.p}
\pl{\sup_{t\in[0,T_u]}Q(t)>u}\sim \pl{Q(0)>u}\sim \pl{\inf_{t\in[0,T_u]}Q(t)>u},
\end{eqnarray}
as $u\to\infty$, providing that $H>1/2$ and $T_u=o(u^{\frac{2H-1}{H}})$; see \cite{PI},
\cite{DE2014}.
Property (\ref{strong.p}) is nowadays referred to as the
{\it strong Piterbarg property}. In \cite{AS2004} it was observed that
(\ref{strong.p}) holds also
for storage processes with self-similar and infinitely divisible input without
Gaussian component.

In this contribution we focus on asymptotic properties of
\BQN
\psi_{T_u}^{\sup}(u):=\pl{\sup_{t\in[0,T_u]}Q(t)>u},
\EQN
and
\BQN
\psi_{T_u}^{\inf}(u):=\pl{\inf_{t\in[0,T_u]}Q(t)>u},
\EQN
as $u\to\infty$,
for wide class of Gaussian processes $X$ and ranges of $T_u$.
As a result, we extend findings of \cite{DI2005}, where
the asymptotics of $\pl{Q(0)>u}$ was considered. Moreover we
generalize \cite{PI} and
\cite{DE2014} where the exact asymptotics of
$\psi_{T_u}^{\sup}(u)$ and $\psi_{T_u}^{\inf}(u)$
were studied for fractional Brownian motion model with $\beta=1$.
As a by-product we find conditions under which
the strong Piterbarg property phenomena (\ref{strong.p}) holds for general  Gaussian $X$
and $\beta$.

Organization of the paper.
Some necessary notation are introduced in Section 2,
whereas the main asymptotic results are presented in Section 3.
In Section 4 we apply derived results to the analysis of
$\psi_{T_u}^{\sup}(u)$ and $\psi_{T_u}^{\inf}(u)$
for $X$ being a sum of independent fractional Brownian motions.
The proofs of main results are given in Section 5. The Appendix contains
proofs of some lemmas that are of technical nature.

\section{Notation}
Throughout this paper we assume that $\{X(t),t\ge0\}$ is
a centered Gaussian process with stationary increments, a.s. continuous sample paths,
$X(0)=0$ and variance function $\sigma^2(t)$ satisfying\\
\\
{\bf AI}: $\sigma^2(t)>0, t>0$ is regularly varying at infinity with index
$2\alpha_\IF\in(0,2)$ and twice continuously differentiable on $(0,\IF)$. Further, its first derivative  $\dot{\sigma^2}$ and second derivative $\ddot{\sigma^2}$ are both ultimately monotone.\\
{\bf AII}: $\sigma^2(t)$ is regularly varying at $0$ with index $2\alpha_0\in(0,2]$.\\
\\
Assumptions {\bf AI-AII} allow us to cover models that play important role in
Gaussian storage models, including both
aggregations of fractional Brownian motions and integrated stationary Gaussian processes;
see, e.g., \cite{Norros94,HP99,DI2005,DE2002}.
{\bf AI-AII} go in line with \cite{DI2005},  where the exact asymptotics of $\pl{Q(0)>u}$, as $u\rw\IF$,
was derived.

Recall that fractional Brownian motion $B_H=\{B_H(t), t\geq 0\}$
with Hurst index  $H\in(0,1)$
is a centered Gaussian process with continuous sample path and covariance function
$
Cov(B_H(t), B_H(s))=\frac{1}{2}\left(|t|^{2H}+|s|^{2H}-|t-s|^{2H}\right).
$

For $X(t)$ satisfying {\bf AI-AII}, $c>0$ and $\beta>\frac{\alpha_\infty}{2}$ we define the storage process
$\{Q(t), t\ge0\}$, where
\BQNY
Q(t)=\sup_{ s \ge t}\left(X(s)-X(t)-c(s-t)^\beta\right),\ \ \ t\geq 0.
\EQNY
Note that assumption $\beta>\frac{\alpha_\infty}{2}$ ensures that $Q(t)$
is finite a.s. for any $t\ge0$.

In order to formulate the main results of this contribution,
following \cite{DE2014}, let
\BQN\label{sv1}
\Phi: C(M)\rw \mathbb{R},
\EQN
be a continuous functional on
the Banach space $C(M)$
of all continuous functions on compact set
$M\subset \mathbb{R}^d, d\geq 1$ with the norm $||f||=\sup_{\vk{t}\in M}|f(\vk{t})|$,
satisfying
\\
{\bf F1:} $|\Phi(f)|\leq \sup_{\vk{t}\in M}f(\vk{t})$.\\
{\bf F2:} $\Phi(af+b)=a\Phi(f)+b,$ for any $a, b>0$.\\
Then,
for a centered continuous Gaussian field
$V=\{V(\vk{t}): \vk{t}\in\mathbb{R}^d\}$ such that
$V(\vk{0})=0$,
\BQN\label{sv}
Cov(V(\vk{t}), V(\vk{s}))=\frac{\sigma_V^2(\vk{t})+\sigma_V^2(\vk{s})-\sigma_V^2(\vk{t}-\vk{s})}{2}
\EQN
and\\
{\bf E1:} $\E{V(\vk{t})-V(\vk{s})}^2=\sigma_V^2(\vk{t}-\vk{s})\leq G|\vk{t}-\vk{s}|^{\alpha_1}$ \\
with $G,\alpha_1>0$,
we introduce
{\it the generalized Pickands' constant}
\BQNY
\mathcal{H}_V^{\Phi}(M)=\E{e^{\Phi\left(\sqrt{2}V-\sigma_{V}^2\right)}}.
\EQNY
We refer to \cite{DE2014} for the finiteness of $\mathcal{H}_V^{\Phi}(M)$.
In particular, for $M=\prod_{i=1}^d[0,S_i]$ with $S_i>0, 1\leq i\leq d$,
and $\Phi(f)=\sup_{\vk{t}\in \prod_{i=1}^d[0,S_i]} f(\vk{t})$,
we use notation
$\mathcal{H}_V(\prod_{i=1}^d[0,S_i])$.
Further, for $d=1$, let
 \BQNY
 \mathcal{H}_V=\lim_{S\rw\IF}\frac{\mathcal{H}_V[0,S]}{S},
 \EQNY
providing that the above limit exists, where $\mathcal{H}_V[0,S]:=\mathcal{H}_V([0,S])$.
We refer to \cite{Pit96}, \cite{DE2002}, \cite{DI2005} and \cite{DE2014}
for the analysis of properties of Pickands'-type constants.

We write
$f_u(t)\RW f(t)$ for $t\in D$ meaning that the convergence is uniform with respect
to $t$ in the domain $D$ as $u\rw\IF$. By $\mathbb{Q}$ and $\mathbb{Q}_i, i=0,1,2,\cdots,n$
we denote some positive constants which may change from line to line.
By $\overleftarrow{\sigma}(\cdot)$ we denote the generalized inverse function to
${\sigma}(\cdot)$,
$\Psi(\cdot)$ denotes the tail distribution of the standard Normal random variable.
We write $f(u)\sim g(u)$ if $\lim_{u\to \infty}\frac{f(u)}{g(u)}=1$.
\def\tA{\Lambda(u)}
\def\tAY{\Lambda^*(u)}

\section{Main Results}\label{s.main}
In this section, we present the exact asymptotics of
$\psi_{T_u}^{\sup}(u)$ and $\psi^{\inf}_{T_u}(u)$.
In further analysis we tacitly assume  that
the variance function $\sigma^2$ of $X$ satisfies both {\bf AI} and {\bf AII}.

Let
$$\varphi:=\lim_{u\rw\IF}\frac{\sigma^2(u^{1/\beta})}{u},$$
assuming that the limit exists.
As it is shown below, according to the value of $\varphi$,
the asymptotics of $\psi_{T_u}^{\sup}(u)$ takes different form.
Additionally, we introduce
$\tau^*=\left(\frac{{\alpha_\IF}}{c(\beta-{\alpha_\IF})}\right)^{1/\beta}$ and
set
\BQN\label{Delta}
\Delta(u):=\left\{
\begin{array}{cc}
\overleftarrow{\sigma}
 \left(\frac{\sqrt{2}\sigma^2(u^{1/\beta}\tau^*)}{u(1+c\tau^{*\beta})}\right),& if \varphi=\IF\ or\  0,\\
 1,&if \varphi\in(0,\IF).\\
 \end{array}
 \right.
\EQN


Let
$$A=\left(\frac{{\alpha_\IF}}{c(\beta-{\alpha_\IF})}\right)^{-\alpha_\IF/\beta}\frac{\beta}{\beta-\alpha_\IF}, \quad
B=\left(\frac{\alpha_\IF}{c(\beta-\alpha_\IF)}\right)^{-(\alpha_\IF+2)/\beta}\alpha_\IF\beta,
\quad m^{*}(u)=\frac{u(1+c(\tau^*)^\beta)}{\sigma(u^{1/\beta})(\tau^*)^{\alpha_\IF}}
.$$

\BT\label{main3} Suppose that $\varphi=0$. \\
i)If  $\frac{T_u}{\Delta(u)}\rw \rho \in [0,\IF)$, then
\BQNY
 \psi_{T_u}^{\sup}(u)\sim \mathcal{H}_{B_{\alpha_0}}[0,\rho]\mathcal{H}_{B_{\alpha_0}}
 \sqrt{\frac{2A\pi}{B}}\frac{u^{1/\beta-1}\sigma(u^{1/\beta}\tau^*)}{(1+c\tau^{*\beta})\Delta(u)}\Psi\left(\inf_{t\geq 0}\frac{u(1+ct^\beta)}{\sigma(u^{1/\beta} t)}\right).
\EQNY
\COM{ii)If $\frac{T_u}{\Delta(u)}\rw \rho\in(0,\IF)$, then
\BQNY
\psi_{T_u}^{\sup}(u)\sim\mathcal{H}_{B_{\alpha_0}}[0,\rho]\mathcal{H}_{B_{\alpha_0}}\sqrt{\frac{2A\pi}{B}}\frac{u^{1/\beta-1}\sigma(u^{1/\beta}\tau^*)}{(1+c\tau^{*\beta})
\overleftarrow{\sigma}
 \left(\frac{\sqrt{2}\sigma^2(u^{1/\beta}\tau^*)}{u(1+c\tau^{*\beta})}\right)}\Psi\left(\inf_{t\geq 0}\frac{u(1+ct^\beta)}{\sigma(u^{1/\beta} t)}\right);
\EQNY
}
ii) If $\frac{T_u}{\Delta(u)}\rw \IF$ and  $T_u=o(e^{\beta_1(m^{*}(u))^2})$  with $\beta_1\in(0, 1/2)$, then
 \BQNY
 \psi_{T_u}^{\sup}(u)\sim\left(\mathcal{H}_{B_{{\alpha_0}}}\right)^2
\sqrt{\frac{2A\pi}{B}}T_u\frac{u^{1/\beta-1}\sigma(u^{1/\beta}\tau^*)}{(1+c\tau^{*\beta})\Delta^2(u)}\Psi\left(\inf_{t\geq 0}\frac{u(1+ct^\beta)}{\sigma(u^{1/\beta} t)}\right) .
 \EQNY
\ET
\BT\label{main4}
 Suppose that $\varphi\in(0,\IF)$.\\
i) If  $T_u\rw \rho\in [0,\IF)$, then
\BQNY
 \psi_{T_u}^{\sup}(u)\sim \mathcal{H}_{\frac{1+c(\tau^*)^\beta}{\sqrt{2}\varphi(\tau^*)^{2\alpha_\IF}}X}[0,\rho] \mathcal{H}_{\frac{1+c(\tau^*)^\beta}{\sqrt{2}\varphi(\tau^*)^{2\alpha_\IF}}X}
 \sqrt{\frac{2A\pi}{B}}\frac{u^{1/\beta-1}\sigma(u^{1/\beta}\tau^*)}{(1+c\tau^{*\beta})\Delta(u)}\Psi\left(\inf_{t\geq 0}\frac{u(1+ct^\beta)}{\sigma(u^{1/\beta} t)}\right).
\EQNY
\COM{ii) If $T_u\rw \rho\in(0,\IF)$, then
\BQNY
\psi_{T_u}^{\sup}(u)\sim \mathcal{H}_{\frac{1+c(\tau^*)^\beta}{\sqrt{2}\varphi(\tau^*)^{2H}}X}[0,\rho] \mathcal{H}_{\frac{1+c(\tau^*)^\beta}{\sqrt{2}\varphi(\tau^*)^{2H}}X}\sqrt{\frac{2A\pi}{B}}\frac{u^{1/\beta-1}\sigma(u^{1/\beta}\tau^*)}{(1+c\tau^{*\beta})
}\Psi\left(\inf_{t\geq 0}\frac{u(1+ct^\beta)}{\sigma(u^{1/\beta} t)}\right);
\EQNY
}
ii) If $T_u\rw \IF$ and  $T_u=o(e^{\beta_1(m^{*}(u))^2})$  with $\beta_1\in(0, 1/2)$, then
 \BQNY
 \psi_{T_u}^{\sup}(u)\sim\left( \mathcal{H}_{\frac{1+c(\tau^*)^\beta}{\sqrt{2}\varphi(\tau^*)^{2\alpha_\IF}}X}\right)^2
 \sqrt{\frac{2A\pi}{B}}T_u\frac{u^{1/\beta-1}\sigma(u^{1/\beta}\tau^*)}{(1+c\tau^{*\beta})\Delta^2(u)}\Psi\left(\inf_{t\geq 0}\frac{u(1+ct^\beta)}{\sigma(u^{1/\beta} t)}\right).
 \EQNY
\ET
\BT\label{main5}  Suppose that $\varphi=\IF$.
\\
i)If  $\frac{T_u}{\Delta(u)}\rw \rho \in [0,\IF)$, then
\BQNY
 \psi_{T_u}^{\sup}(u)\sim \mathcal{H}_{B_{\alpha_\infty}}[0,\rho]\mathcal{H}_{B_{\alpha_\infty}}
 \sqrt{\frac{2A\pi}{B}}\frac{u^{1/\beta-1}\sigma(u^{1/\beta}\tau^*)}{(1+c\tau^{*\beta})\Delta(u)}\Psi\left(\inf_{t\geq 0}\frac{u(1+ct^\beta)}{\sigma(u^{1/\beta} t)}\right).
\EQNY
\COM{ii)If $\frac{T_u}{\Delta(u)}\rw \rho\in(0,\IF)$, then
\BQNY
\psi_{T_u}^{\sup}(u)\sim\mathcal{H}_{B_{\alpha_0}}[0,\rho]\mathcal{H}_{B_{\alpha_0}}\sqrt{\frac{2A\pi}{B}}\frac{u^{1/\beta-1}\sigma(u^{1/\beta}\tau^*)}{(1+c\tau^{*\beta})
\overleftarrow{\sigma}
 \left(\frac{\sqrt{2}\sigma^2(u^{1/\beta}\tau^*)}{u(1+c\tau^{*\beta})}\right)}\Psi\left(\inf_{t\geq 0}\frac{u(1+ct^\beta)}{\sigma(u^{1/\beta} t)}\right);
\EQNY
}
ii) If $\frac{T_u}{\Delta(u)}\rw \IF$ and  $T_u=o(e^{\beta_1(m^{*}(u))^2})$  with $\beta_1\in(0, 1/2)$, then
 \BQNY
 \psi_{T_u}^{\sup}(u)\sim\left(\mathcal{H}_{B_{{\alpha_\infty}}}\right)^2
 \sqrt{\frac{2A\pi}{B}}T_u\frac{u^{1/\beta-1}\sigma(u^{1/\beta}\tau^*)}{(1+c\tau^{*\beta})\Delta^2(u)}\Psi\left(\inf_{t\geq 0}\frac{u(1+ct^\beta)}{\sigma(u^{1/\beta} t)}\right) .
 \COM{\sqrt{\frac{2A\pi}{B}}\frac{T_uu^{1/\beta-1}\sigma(u^{1/\beta}\tau^*)}{(1+c\tau^{*\beta})
 \left(\overleftarrow{\sigma}
 \left(\frac{\sqrt{2}\sigma^2(u^{1/\beta}\tau^*)}{u(1+c\tau^{*\beta})}\right)\right)^2}\Psi\left(\inf_{t\geq 0}\frac{u(1+ct^\beta)}{\sigma(u^{1/\beta} t)}\right).
 }
 \EQNY
\ET
The above trichotomy with respect to the value of $\varphi$ goes in line with
findings of Dieker \cite{DI2005}, where the asymptotics of $\pl{Q(0)>u}$, as $u\to\infty$,
was derived.

The following theorem deals with the asymptotic behavior of the tail distribution
of $\psi^{\inf}_{T_u}(u)$.

\BT\label{main6} i) If $\varphi=0$  and $\frac{T_u}{\Delta(u)}\rw \rho \in [0,\IF)$, then
\BQNY
 \psi^{\inf}_{T_u}(u)\sim \mathcal{H}^{\inf}_{B_{\alpha_0}}[0,\rho]\mathcal{H}_{B_{\alpha_0}}
  \sqrt{\frac{2A\pi}{B}}\frac{u^{1/\beta-1}\sigma(u^{1/\beta}\tau^*)}{(1+c\tau^{*\beta})\Delta(u)}\Psi\left(\inf_{t\geq 0}\frac{u(1+ct^\beta)}{\sigma(u^{1/\beta} t)}\right);
\EQNY
ii) If $\varphi\in(0,\IF)$  and $T_u\rw \rho \in[0,\IF)$, then
\BQNY
 \psi^{\inf}_{T_u}(u)\sim \mathcal{H}^{\inf}_{\frac{1+c(\tau^*)^\beta}{\sqrt{2}\varphi(\tau^*)^{2\alpha_\IF}}X}[0,\rho] \mathcal{H}_{\frac{1+c(\tau^*)^\beta}{\sqrt{2}\varphi(\tau^*)^{2\alpha_\IF}}X}
  \sqrt{\frac{2A\pi}{B}}\frac{u^{1/\beta-1}\sigma(u^{1/\beta}\tau^*)}{(1+c\tau^{*\beta})\Delta(u)}\Psi\left(\inf_{t\geq 0}\frac{u(1+ct^\beta)}{\sigma(u^{1/\beta} t)}\right);
 \EQNY
iii) If $\varphi=\IF$ and $\frac{T_u}{\Delta(u)}\rw \rho\in [0,\IF)$, then
\BQNY
 \psi^{\inf}_{T_u}(u)\sim \mathcal{H}^{\inf}_{B_{\alpha_\IF}}[0,\rho]\mathcal{H}_{B_{\alpha_\IF}}
  \sqrt{\frac{2A\pi}{B}}\frac{u^{1/\beta-1}\sigma(u^{1/\beta}\tau^*)}{(1+c\tau^{*\beta})\Delta(u)}\Psi\left(\inf_{t\geq 0}\frac{u(1+ct^\beta)}{\sigma(u^{1/\beta} t)}\right).
\EQNY
\ET
Combination of the above findings straightforwardly leads to the following corollary that
deals with the {\it strong Piterbarg property} for $Q$, extending results derived in \cite{DE2014}.
\BK \label{cor.1}
Suppose that $\frac{T_u}{\Delta(u)}\rw0$.

Then
$$\psi_{T_u}^{\inf}(u)\sim\psi_{T_u}^{\sup}(u)\sim \psi_{0}(u).$$
\EK

\begin{remark}
The relation $\frac{T_u}{\Delta(u)}\rw0$ in Corollary \ref{cor.1}
is optimal. Indeed, if $\frac{T_u}{\Delta(u)}\rw \rho>0$, then comparing Theorems
\ref{main3}, \ref{main4}, \ref{main5} and \ref{main6},
none of the asymptotic relation in Corollary \ref{cor.1} holds.
\end{remark}
\section{Application to heterogenous fluid queues}
Consider the stationary storage model
\BQN
Q(t)=\sup_{s\geq t}\left(\sum_{i=1}^n\left(B_{H_i}(s)-B_{H_i}(t)\right)-c(s-t)^\beta\right),\ \ \ t\geq 0,
\EQN
where
$B_{H_i}(t), 1\leq i\leq n$ are mutually independent fractional Brownian motions with indexes
$1>H_1>H_2\geq\cdots\geq H_{n-1}>H_n>0$ respectively and $\beta>H_1$.
It is straightforward to check that
$\sigma^2_\Sigma(t):=Var\left(\sum_{i=1}^nB_{H_i}(t)\right)=\sum_{i=1}^nt^{2H_i}$
satisfies {\bf AI-AII} with $\alpha_0=2H_n$ and $\alpha_\infty=2H_1$,
which in the light of Theorems \ref{main3}, \ref{main4} and \ref{main5}, leads to.
\BK Suppose that $2H_1<\beta$.\\
\COM{i)If  $T_u=o(u^{\frac{2H_1-\beta}{\beta H_n}})$, then
\BQNY
 \psi_{T_u}^{\sup}(u)\sim \mathcal{H}_{B_{H_n}}2^{\frac{H_n-1}{2H_n}}\sqrt{\frac{A\pi}{B}}\frac{(1+c(\tau^*)^\beta)^{\frac{1-H_n}{H_n}}}{(\tau^*)^{\frac{H_1(2-H_n)}{H_n}}}
 u^{\frac{H_n+\beta-2H_1+H_1H_n-\beta H_n}{\beta H_n}}\Psi\left(\inf_{t\geq 0}\frac{u(1+ct^\beta)}{\sigma_\Sigma(u^{1/\beta} t)}\right);
\EQNY}
i)If $T_u u^{\frac{\beta-2H_1}{\beta H_n}}
\sim \rho \left(\frac{\sqrt{2}(\tau^*)^{2H_1}}{1+c(\tau^*)^\beta}\right)^{\frac{1}{H_n}}$ as $u\to\infty$, with $\rho\in [0,\IF)$, then
\BQNY
\psi_{T_u}^{\sup}(u)\sim\mathcal{H}_{B_{H_n}}[0,\rho]\mathcal{H}_{B_{H_n}}2^{\frac{H_n-1}{2H_n}}\sqrt{\frac{A\pi}{B}}\frac{(1+c(\tau^*)^\beta)^{\frac{1-H_n}{H_n}}}
{(\tau^*)^{\frac{H_1(2-H_n)}{H_n}}}
 u^{\frac{H_n+\beta-2H_1+H_1H_n-\beta H_n}{\beta H_n}}\Psi\left(\inf_{t\geq 0}\frac{u(1+ct^\beta)}{\sigma_\Sigma(u^{1/\beta} t)}\right);
\EQNY
ii) If $T_u u^{\frac{\beta-2H_1}{\beta H_n}}\rw \IF$ and  $T_u=o(e^{\beta_1(m^{*}(u))^2})$  with $\beta_1\in(0, 1/2)$, then
 \BQNY
 \psi_{T_u}^{\sup}(u)\sim(\mathcal{H}_{B_{H_n}})^22^{\frac{H_n-2}{2H_n}}\sqrt{\frac{A\pi}{B}}\frac{(1+c(\tau^*)^\beta)^{\frac{2-H_n}{H_n}}}
{(\tau^*)^{\frac{H_1(4-H_n)}{H_n}}}
T_u u^{\frac{H_n+2\beta-4H_1+H_1H_n-\beta H_n}{\beta H_n}}\Psi\left(\inf_{t\geq 0}\frac{u(1+ct^\beta)}{\sigma_\Sigma(u^{1/\beta} t)}\right).
 \EQNY
\EK
\BK
Suppose that $2H_1=\beta$.\\
\COM{i) If  $T_u\rw 0$, then
\BQNY
 \psi_{T_u}^{\sup}(u)\sim \mathcal{H}_{\frac{1+c(\tau^*)^\beta}{\sqrt{2}(\tau^*)^{2H}}X_F}\sqrt{\frac{2A\pi}{B}}\frac{(\tau^*)^{H_1}}{(1+c\tau^{*\beta})}u^{\frac{1-\beta+H_1}{\beta}}
 \Psi\left(\inf_{t\geq 0}\frac{u(1+ct^\beta)}{\sigma_\Sigma(u^{1/\beta} t)}\right);
\EQNY}
i) If $T_u\rw \rho\in[0,\IF)$, then
\BQNY
\psi_{T_u}^{\sup}(u)\sim \mathcal{H}_{\frac{1+c(\tau^*)^\beta}{\sqrt{2}(\tau^*)^{2H}}\sum_{i=1}^n B_{H_i}}[0,\rho] \mathcal{H}_{\frac{1+c(\tau^*)^\beta}{\sqrt{2}(\tau^*)^{2H}}\sum_{i=1}^n B_{H_i}}\sqrt{\frac{2A\pi}{B}}\frac{(\tau^*)^{H_1}}{(1+c\tau^{*\beta})}u^{\frac{1-H_1}{2H_1}}
\Psi\left(\inf_{t\geq 0}\frac{u(1+ct^\beta)}{\sigma_\Sigma(u^{1/\beta} t)}\right);
\EQNY
ii) If $T_u\rw \IF$ and  $T_u=o(e^{\beta_1(m^{*}(u))^2})$  with $\beta_1\in(0, 1/2)$, then
 \BQNY
 \psi_{T_u}^{\sup}(u)\sim\left( \mathcal{H}_{\frac{1+c(\tau^*)^\beta}{\sqrt{2}(\tau^*)^{2H}}\sum_{i=1}^n B_{H_i}}\right)^2\sqrt{\frac{2A\pi}{B}}\frac{(\tau^*)^{H_1}}
 {(1+c\tau^{*\beta})}T_uu^{\frac{1-H_1}{2H_1}}
 \Psi\left(\inf_{t\geq 0}\frac{u(1+ct^\beta)}{\sigma_\Sigma(u^{1/\beta} t)}\right).
 \EQNY
\EK
\BK
Suppose that $2H_1>\beta>H_1$.\\
\COM{i) If  $T_u u^{\frac{\beta-2H_1}{\beta H_1}}\rw 0$, then
\BQNY
 \psi_{T_u}^{\sup}(u)\sim\mathcal{H}_{B_{H_1}}2^{\frac{H_1-1}{2H_1}}\sqrt{\frac{A\pi}{B}}\frac{(1+c(\tau^*)^\beta)^{\frac{1-H_1}{H_1}}}{(\tau^*)^{2-H_1}}
 u^{\frac{(\beta-H_1)(1-H_1)}{\beta H_1}}
 \Psi\left(\inf_{t\geq 0}\frac{u(1+ct^\beta)}{\sigma_\Sigma(u^{1/\beta} t)}\right);
\EQNY}
i) If $T_u u^{\frac{\beta-2H_1}{\beta H_1}}\rw \rho \left(\frac{\sqrt{2}(\tau^*)^{2H_1}}{1+c(\tau^*)^\beta}\right)^{\frac{1}{H_1}}$, with $\rho\in[0,\IF)$, then
\BQNY
\psi_{T_u}^{\sup}(u)\sim\mathcal{H}_{B_{H_1}}[0,\rho] \mathcal{H}_{B_{H_1}}2^{\frac{H_1-1}{2H_1}}\sqrt{\frac{A\pi}{B}}\frac{(1+c(\tau^*)^\beta)^{\frac{1-H_1}{H_1}}}{(\tau^*)^{2-H_1}}
 u^{\frac{(\beta-H_1)(1-H_1)}{\beta H_1}}
 \Psi\left(\inf_{t\geq 0}\frac{u(1+ct^\beta)}{\sigma_\Sigma(u^{1/\beta} t)}\right);
\EQNY
ii) If $T_u u^{\frac{\beta-2H_1}{\beta H_1}}\rw \IF$ and  $T_u=o(e^{\beta_1(m^{*}(u))^2})$  with $\beta_1\in(0, 1/2)$, then
 \BQNY
 \psi_{T_u}^{\sup}(u)\sim\left( \mathcal{H}_{B_{H_1}}\right)^22^{\frac{H_1-2}{2H_1}}\sqrt{\frac{A\pi}{B}}\frac{(1+c(\tau^*)^\beta)^{\frac{2-H_1}{H_1}}}{(\tau^*)^{4-H_1}}
 u^{\frac{2\beta-3H_1+H_1^2-\beta H_1}{\beta H_1}}
 \Psi\left(\inf_{t\geq 0}\frac{u(1+ct^\beta)}{\sigma_\Sigma(u^{1/\beta} t)}\right).
 \EQNY
\EK
\BRM
Following  \cite{DE2014} and \cite{PI}, if
$n=1$, $H_1>1/2$, $\beta=1$ and $T_u=o(u^{\frac{2H_1-1}{H_1}})$, then
$ \psi_{T_u}^{\sup}(u)\sim \psi^{\inf}_{T_u}(u)$.
Combination of results derived in Section \ref{s.main} to the model considered in this section
extends this findings to $n\ge1$, $2H_1>\beta>H_1$ and $T_u=o(u^{\frac{2H_1-\beta}{\beta H_1}})$.
\ERM
\section{Proofs}
In this section we present detailed proofs of the main
results of this contribution.

Following the same line of reasoning as in \cite{PI}, we write
\BQNY
\psi_{T_u}^{\sup}(u)=\pl{\sup_{t\in[0,T_u]}Q (t)>u}=\pl{\sup_{t\in[0,u^{-1/\beta}T_u]}\sup_{s\geq t}Z_u(s,t)>m(u)}
\EQNY
with $Z_u(s,t)=\frac{X(u^{1/\beta}s)-X(u^{1/\beta}t)}
{1+c(s-t)^\beta } \frac{1+c\tau_u^\beta}{\sigma(u^{1/\beta}\tau_u)}$
and
$m(u)=\inf_{t\geq 0}\frac{u(1+ct^\beta)}{\sigma(u^{1/\beta}t)}
$.

Hereafter, for a given process $Y(t)$, we denote $\overline{Y}(t):=Y(t)/ \sigma_Y(t)$.
By $\dot{h}$, $\ddot{h}$ we  mean the first and second derivative of twice continuously differentiable function $h$,  respectively.
To short the notation we set
$\sigma_u^2(s)=\E{\frac{X(u^{1/\beta}s)}{\sigma(u)(1+cs^\beta)}}^2$

and $r_u(s,t,s_1,t_1):= Cov(\overline{Z_u}(s,t), \overline{Z_u}(s_1,t_1))=
\E{\frac{X(u^{1/\beta}s)-X(u^{1/\beta}t)}{\sigma(u^{1/\beta}(s-t))} \frac{X(u^{1/\beta}s_1)-X(u^{1/\beta}t_1)}{\sigma(u^{1/\beta}(s_1-t_1))}}, s>t, s_1>t_1$.

The following lemma slightly extends
Lemma 2 in \cite{DI2005}, by providing asymptotics for
the tail distribution of functionals introduced in (\ref{sv1}) and fulfilling {\bf F1-F2}
instead of $\sup$ functional considered in \cite{DI2005}.
Following the
setting given in \cite{DI2005}, let
$\{K_u\}$ be a nondecreasing family of subsets
of $\mathbb{Z}^m$ with $m\geq 1$,
and $\{X^{(u,\vk{k})}(\vk{t}), \vk{t}\in M\}, u>0, \vk{k}\in K_u$ be a collection of
centered continuous Gaussian fields on a compact set $M\subset \mathbb{R}^d$. We assume that the variance of
$X^{(u,\vk{k})}(\vk{t})$ equals 1. Let $g_{\vk{k}},\theta_{\vk{k}}$, with $\vk{k}\in K_u$ be such that (see \cite{DI2005})\\
{\bf P1} $\inf_{\vk{k}\in K_u}g_{\vk{k}}(u)\rw\IF$ as $u\rw\IF$.\\
{\bf P2} There exists a centered Gaussian field $\{V(\vk{t}), {\vk{t}}\in \mathbb{R}^d\}$ with covariance as in (\ref{sv}), satisfying {\bf E1}, such that
$\sup_{\vk{k}\in K_u}|\theta_{\vk{k}}(u,\vk{s},\vk{t})-\sigma_V^2(\vk{t}-\vk{s})|\rw 0$ for any $\vk{s},\vk{t}\in M$.\\
{\bf P3} For some $\eta_1, \cdots, \eta_d>0$,
\BQNY
\limsup_{u\rw\IF}\sup_{\vk{k}\in K_u}\sup_{\vk{s},\vk{t}\in M}\frac{\theta_{\vk{k}}(u,\vk{s},\vk{t})}{\sum_{i=1}^d|s_i-t_i|^{\eta_i}}<\IF.
\EQNY
{\bf P4}
\BQNY
\lim_{\epsilon\rw 0}\limsup_{u\rw\IF}\sup_{\vk{k}\in K_u}\sup_{|\vk{s}-\vk{t}|<\epsilon, \vk{s},\vk{t}\in M}g_{\vk{k}}^2(u)\E{\left(X^{(u,\vk{k})}(\vk{s})-X^{(u,\vk{k})}(\vk{t})\right)X^{(u,\vk{k})}(\vk{0})}=0.
\EQNY
\BEL\label{PICKANDS}
Suppose that
{\bf P1-P4} hold for functions $g_{\vk{k}}$, $\theta_{\vk{k}}$ and Gaussian process $V$.
Let
$\Phi:C(M)\rw \mathbb{R}$ be a continuous functional fulfilling {\bf F1} and {\bf F2}. If
\BQN\label{Cor}
\lim_{u\rw\IF}\sup_{\vk{k}\in K_u}\sup_{\vk{s},\vk{t}\in M, \vk{s}\neq \vk{t}}\left|g_{\vk{k}}^2(u)\frac{Var\left(X^{(u,\vk{k})}(\vk{t})-X^{(u,\vk{k})}(\vk{s})\right)}{2\theta_{\vk{k}}(u,\vk{s},\vk{t})}-1\right| = 0,
\EQN
then 
$$\limit{u} \frac{\mathbb{P}\left(\Phi(X^{(u,\vk{k})})>g_{\vk{k}}(u)\right)}{\Psi(g_{\vk{k}}(u))}=\mathcal{H}_{V}^{\Phi}(M)
$$
provided that $\mathbb{P}\left(\Phi(X^{(u,\vk{k})})>g_{\vk{k}}(u)\right)>0$ for $u$ large enough,
and
\begin{eqnarray}\label{PB}
\limsup_{u\to\infty} \sup_{\vk{k} \in K_u} \frac{\mathbb{P}\left(\Phi(X^{(u,\vk{k})})>g_{\vk{k}}(u)\right)}{\Psi(g_{\vk{k}}(u))}<\infty.
\end{eqnarray}
\EEL
The proof of Lemma \ref{PICKANDS}
goes line-by line the same as the proof of Lemma 2 in \cite{DI2005}; see also proof of Lemma 1 in \cite{DE2014}. We
present main steps of the proof in Appendix.

\BEL\label{R1}
Suppose that $\sigma^2(t)$ satisfies {\bf AI-AII}. Then there exisit $\gamma\in (0,2)$, $C_1>0$ and $C_2>0$ such that
 \BQNY
 C_1t^2\leq\sigma^2(t)\leq C_2t^\gamma,
 \EQNY
 holds in a neighbourhood of zero.
\EEL

\BEL\label{L1}
 For $u$ large enough,   $\sigma_u(\tau)$ attains its unique maximum $\tau_u \in [0,\IF)$
 so that $\tau_u\rw \tau^*$, as
 $u\to \infty$. Moreover,
\BQN\label{V4}
\frac{\sigma_u(\tau)}{\sigma_u(\tau_u)}=1-b_u(\tau-\tau_u)^2(1+o(1)),\ \ \ \tau\rw\tau_u,
\EQN
where $b_u\rw b=\frac{B}{2A}$.
\EEL
Let $E(u)=(\tau_u-\delta_u,\tau_u+\delta_u)$ with
$\delta_u=\frac{\ln (m(u))}{m(u)}$.

\BEL\label{L2}
 The correlation function $r_u(s,t,s_1,t_1)$ satisfies
\BQNY
\lim_{u\rw \IF}\sup_{|t-t_1|<\delta_u,s-t,s_1-t_1\in E(u),(s,t)\neq (s_1,t_1)}\left|\frac{1-r_u(s,t,s_1,t_1)}{\frac{\sigma^2(u^{1/\beta}|s-s_1|)+\sigma^2(u^{1/\beta}|t-t_1|)}{2\sigma^2(u^{1/\beta}\tau^*)}}-1\right|=0.
\EQNY
\EEL

\BEL\label{L3}
For $u$ large enough and any $\delta>0$, there exists a constant $0<a_\delta<1$ such that $$r_u(s,t,s_1,t_1)<a_\delta$$ holds for all $|t-t_1|>\delta ,s-t,s_1-t_1\in E(u)$.
Further, $$\lim_{R\rw\IF}\sup_{|t-t_1|>R, s-t,s_1-t_1\in E(u)}r_u(s,t,s_1,t_1)=0,$$
holds uniformly with respect to $u$ for $u$ large enough.
\EEL
We provide complete proofs of Lemma \ref{PICKANDS}, \ref{R1}, \ref{L1}, \ref{L2} and \ref{L3} in the Appendix.

The following lemma deals with the asymptotics of supremum of Gaussian field $Z_u(s,t)$
over a parameter set that is away of the neighbourhoood of the maximizer of the variance of  $Z_u(s,t)$.

Recall that
\BQNY
\Delta(u)=\left\{
\begin{array}{cc}
\overleftarrow{\sigma}
 \left(\frac{\sqrt{2}\sigma^2(u^{1/\beta}\tau^*)}{u(1+c\tau^{*\beta})}\right),& if \varphi=\IF\ or\  0,\\
 1,&if \varphi\in(0,\IF).\\
 \end{array}
 \right.
\EQNY
\BEL\label{L4} Suppose that {\bf AI-\bf AII} hold and  $r>0$.\\
i) If $T_u=o(u^r)$, then
\BQNY
\pl{\sup_{t\in[0,T_uu^{-1/\beta}]}\sup_{s\geq t, s-t\notin E(u)}Z_u(s,t)>m(u)}=o\left(\frac{u^{1/\beta}}{\Delta(u)m(u)}\Psi(m(u))\right).
\EQNY
ii) If $u^rT_u\rw \IF$, as $u\rw \IF$, then
\BQNY
\pl{\sup_{t\in[0,T_uu^{-1/\beta}]}\sup_{s\geq t, s-t\notin E(u)}Z_u(s,t)>m(u)}=o\left(\frac{T_uu^{1/\beta}}{(\Delta(u))^2m(u)}\Psi(m(u))\right).
\EQNY
\EEL
{\bf Proof.}
We set $\tau=s-t$ and write
\[
\pl{\sup_{t\in[0,T_uu^{-1/\beta}]}\sup_{s-t\notin E(u)}Z_u(s,t)>m(u)}
=
\pl{\sup_{t\in[0,T_uu^{-1/\beta}]}\sup_{\tau\notin E(u)}Z_u(t+\tau,t)>m(u)}.
\]
Let
$
[0,T_uu^{-1/\beta}]\times ([0,\infty)\setminus E(u))=
\mathcal{S}_{1,u}\cup \mathcal{S}_{2,u}\cup\mathcal{S}_{3,u},
$
where
$\mathcal{S}_{1,u}=[0,T_uu^{-1/\beta}]\times[0,\epsilon]$,
$\mathcal{S}_{2,u}=[0,T_uu^{-1/\beta}]\times[T,\IF)$
and
$\mathcal{S}_{3,u}=[0,T_uu^{-1/\beta}]\times([\epsilon,T]\setminus E(u))$, for sufficiently small $\epsilon>0$ and $T\in \mathbb{N}^+$.
Clearly it suffices to find the asymptotic upper estimates of the analyzed tail probability
for each set $\mathcal{S}_{1,u},\mathcal{S}_{2,u},\mathcal{S}_{3,u}$ separately.

{\it Ad.  $ \mathcal{S}_{2,u}$.}
Following Potter's theorem (see, e.g., \cite{BI1989}),
\BQN\label{V8}
\E{Z_u^2(t+\tau,t)}&=&\frac{\sigma^2(u^{1/\beta}\tau)(1+c\tau_u^\beta)^2}{\sigma^2(u^{1/\beta}\tau_u)(1+c\tau^\beta)^2}\leq 2\left(\frac{\tau}{\tau_u}\right)^{2{{\alpha_\IF}}+2\eta}\left(\frac{1+c\tau^\beta_u}{1+c\tau^\beta}\right)^2\leq 2\frac{(1+c\tau^\beta_u)^2}{c^2\tau_u^{2{\alpha_\IF}+2\eta}}\tau^{2{\alpha_\IF}+2\eta-2\beta}\nonumber\\
&\leq& \mathbb{Q}\tau^{-2(\beta-{\alpha_\IF}-\eta)},
\EQN
where $\mathbb{Q}$ is a fixed constant and $0<\eta<\beta-{\alpha_\IF}$.
Note that by {\bf AI} and Remark \ref{R1}, we can choose $0<\gamma_1<\min(\gamma,2{\alpha_\IF})$ such that
\BQN\label{EG}
g_1(t):=\frac{\sigma^2(t)}{t^{\gamma_1}}
\EQN
is a regularly varying function at $\IF$ with index $2{\alpha_\IF}-\gamma_1>0$ and bounded in a neighborhood of $0$. Therefore it follows from
Uniform Convergence Theorem (UCT) (see, e.g., Theorem 1.5.2 in \cite{BI1989})
that for $t,t_1\in [l,l+1]\subset[0,T_uu^{-1/\beta}+1]$ and $\tau,\tau_1\in[k,k+1]\subset[T,\IF)$,
\BQN\label{V9}
\lefteqn{\E{\overline{Z}_u(t+\tau,t)-\overline{Z}_u(t_1+\tau_1,t_1)}^2}\nonumber\\
&\leq& 2\frac{\sigma^2(u^{1/\beta}|t-t_1|)+\sigma^2(u^{1/\beta}|t+\tau-t_1-\tau_1|)}{\sigma(u^{1/\beta}\tau)\sigma(u^{1/\beta}\tau_1)}\nonumber\\
&\leq& 4\frac{g_1(u^{1/\beta}|t-t_1|)}{g_1(u^{1/\beta}k)}\left(\frac{|t-t_1|}{k}\right)^{\gamma_1}+4\frac{g_1(u^{1/\beta}|t+\tau-t_1-\tau_1|)}
{g_1(u^{1/\beta}k)}\left(\frac{|t+\tau-t_1-\tau_1|}{k}\right)^{\gamma_1}\nonumber\\
&\leq&\mathbb{Q}_1(|t-t_1|^{\gamma_1}+|t+\tau-t_1-\tau_1|^{\gamma_1})
\leq\mathbb{Q}_2(|t-t_1|^{\gamma_1}+|\tau-\tau_1|^{\gamma_1}),
\EQN
where $\mathbb{Q}_1,\mathbb{Q}_2>0$.
Combining (\ref{V8}) and (\ref{V9}) with Fernique inequality (see \cite{LE1983}), we have
\BQN\label{E3}
\lefteqn{\pl{\sup_{t\in[0,T_uu^{-1/\beta}]}\sup_{\tau\in[T,\IF)}Z_u(t+\tau,t)>m(u)}}\nonumber\\
&\leq& \sum_{l=0}^{[T_uu^{-1/\beta}]}\sum_{k=T}^\IF\pl{\sup_{t\in[l,l+1]}\sup_{\tau\in[k,k+1] }Z_u(t+\tau,t)>m(u)}\nonumber\\
&\leq&\sum_{l=0}^{[T_uu^{-1/\beta}]}\sum_{k=T}^\IF\pl{\sup_{t\in[l,l+1]}\sup_{\tau\in[k,k+1] }\overline{Z}_u(t+\tau,t)>m(u)\frac{k^{(\beta-{\alpha_\IF}-\eta)}}{\sqrt{\mathbb{Q}}}}\nonumber\\
&\leq&\sum_{l=0}^{[T_uu^{-1/\beta}]}\sum_{k=T}^\IF 16\exp\left(-\frac{k^{2(\beta-{\alpha_\IF}-\eta)}}{\mathbb{Q}_3}m^2(u)\right)\nonumber\\
&\leq& \mathbb{Q}_T([T_uu^{-1/\beta}]+1)\exp\left(-\frac{T^{2(\beta-{\alpha_\IF}-\eta)}}{\mathbb{Q}_3}m^2(u)\right),
\EQN
for $T$ large enough, with $\mathbb{Q}_3,\mathbb{Q}_T>0$.\\
{\it Ad.  $ \mathcal{S}_{1,u}$.}
It is convenient to bound
\BQN
\pl{\sup_{t\in[0,T_uu^{-1/\beta}],\tau\in[0,\epsilon]}Z_u(t+\tau,t)>m(u)}&\leq& \pl{\sup_{t\in[0,T_uu^{-1/\beta}],\tau\in[0,\epsilon]}Z_u(t+\tau,t)(1+c\tau^\beta)>m(u)}.
\EQN
Indeed the same lines of reasoning as above leads 
for $t,t_1\in [l,l+1]\subset [0,T_uu^{-1/\beta}+1]$ and $\tau, \tau_1\in[0,\epsilon]$, to
\BQNY
\lefteqn{\E{Z_u(t+\tau,t)(1+c\tau^\beta)-Z_u(t_1+\tau_1,t_1)(1+c\tau_1^\beta)}^2}\\
&=&\frac{(1+c\tau_u^\beta)^2}{\sigma^2(u^{1/\beta}\tau_u)}
\E{X(u^{1/\beta}(t+\tau))-X(u^{1/\beta}t)-X(u^{1/\beta}(t_1+\tau_1))+X(u^{1/\beta}t_1)}^2\\
&\leq& 2(1+c\tau_u^\beta)^2\frac{\sigma^2(u^{1/\beta}|t-t_1|)+\sigma^2(u^{1/\beta}|t+\tau-t_1-\tau_1|)}{\sigma^2(u^{1/\beta}\tau_u)}\\
&=&2(1+c\tau_u^\beta)^2\left(\frac{g_1(u^{1/\beta}|t-t_1|)}{g_1(u^{1/\beta}\tau_u)}\frac{|t-t_1|^{\gamma_1}}{\tau_u^{\gamma_1}}
+\frac{g_1(u^{1/\beta}|t+\tau-t_1-\tau_1|)}{g_1(u^{1/\beta}\tau_u)}\frac{|t+\tau-t_1-\tau_1|^{\gamma_1}}{\tau_u^{\gamma_1}}\right)\\
&\leq& \mathbb{Q}_4 \left(|t-t_1|^{\gamma_1}+|\tau-\tau_1|^{\gamma_1}\right),
\EQNY
where $\mathbb{Q}_4>0$.
Therefore, by Fernique inequality,
$\pl{\sup_{t\in[l,l+1],\tau\in[0,\epsilon]}Z_u(t+\tau,t)(1+c\tau^\beta)>x}\leq 8\exp{\left(-\frac{x^2}{\mathbb{Q}_5}\right)},
$
for all $[l,l+1]\subset[0,T_uu^{-1/\beta}+1]$ and $x>0$.
Hence we can find a common $a>0$ such that
\\
$\pl{\sup_{t\in[l,l+1],\tau\in[0,\epsilon]}Z_u(t+\tau,t)(1+c\tau^\beta)>a}< 1/2,
$
for all $[l,l+1]\subset[0,T_uu^{-1/\beta}+1]$. Moreover, we have
\BQNY
\sup_{\tau\in[0,\epsilon]}\E{Z_u(t+\tau,t)(1+c\tau^\beta)}^2\leq \sup_{\tau\in[0,\epsilon]}\frac{\sigma^2(u^{1/\beta}\tau)(1+c\tau_u^{\beta})^2}{\sigma^2(u^{1/\beta}\tau_u)}
\leq\mathbb{Q}_6\left(\frac{\epsilon}{\tau^*}\right)^{2{\alpha_\IF}}.
\EQNY
Thus, by Borel theorem, we have that for $\epsilon$ small enough,
\BQN\label{E4}
\lefteqn{\pl{\sup_{t\in[0,T_uu^{-1/\beta}],\tau\in[0,\epsilon]}Z_u(t+\tau,t)(1+c\tau^\beta)>m(u)}}\nonumber\\
&\leq&\sum_{l=0}^{[T_uu^{-1/\beta}]}\pl{\sup_{t\in[l,l+1],\tau\in[0,\epsilon]}Z_u(t+\tau,t)(1+c\tau^\beta)>m(u)}\nonumber\\
&\leq&2([T_uu^{-1/\beta}]+1)\Psi\left(\frac{m(u)-a}{\sqrt{\mathbb{Q}_6}\left(\frac{\epsilon}{\tau^*}\right)^{{\alpha_\IF}}}\right).
\EQN
{\it Ad.  $ \mathcal{S}_{3,u}$.}
Similarly as for $\mathcal{S}_{2,u}$, we have
\BQNY
\E{\overline{Z}_u(t+\tau,t)-\overline{Z}_u(t_1+\tau_1,t_1)}^2\leq \mathbb{Q}_7(|t-t_1|^{\gamma_1}+|\tau-\tau_1|^{\gamma_1}),
\EQNY
for $t,t_1\in [l,l+1]\subset [0,T_uu^{-1/\beta}+1]$ and $\tau, \tau_1\in[\epsilon,T]$.
Thus by Piterbarg inequality (Theorem 8.1 in \cite{Pit96}) and (\ref{V4}), for any $\epsilon>0$, we have
\BQN\label{E5}
\lefteqn{\pl{\sup_{(t,\tau)\in[0,T_uu^{-1/\beta}]\times([\epsilon,T]/ E(u))}Z_u(t+\tau,t)>m(u)}}\nonumber\\
&\leq& \pl{\sup_{(t,\tau)\in[0,T_uu^{-1/\beta}]\times[\epsilon,T]}\overline{Z}_u(t+\tau,t)>\frac{m(u)}{1-\frac{b}{2}\left(\frac{\ln m(u)}{m(u)}\right)^2}}\nonumber\\
&\leq&\mathbb{Q}_8([T_uu^{-1/\beta}]+1)T(m(u))^{4/\gamma_1}\Psi\left(\frac{m(u)}{1-\frac{b}{2}\left(\frac{\ln m(u)}{m(u)}\right)^2}\right)=o\left(\frac{u^{1/\beta-\epsilon}}{\Delta(u)m(u)}\Psi(m(u))\right).
\EQN
Combination of  (\ref{E3}), (\ref{E4}) and (\ref{E5}) establishes the claims.\QED\\

\subsection{Proof of Theorem 3.1}
Since
\BQN\label{F0}
\pi_{T_u}(u)\leq \psi_{T_u}^{\sup}(u)\leq \pi_{T_u}(u)+\pl{\sup_{t\in[0,T_uu^{-1/\beta}]}\sup_{s\geq t, s-t\notin E(u)}Z_u(s,t)>m(u)},
\EQN
where
\BQNY
\pi_{T_u}(u)=\pl{\sup_{t\in[0,T_uu^{-1/\beta}]}\sup_{s-t\in E(u)}Z_u(s,t)>m(u)}.
\EQNY
and the upper bound of
$\pl{\sup_{t\in[0,T_uu^{-1/\beta}]}\sup_{s\geq t, s-t\notin E(u)}Z_u(s,t)>m(u)}$
is given in Lemma {\ref{L4}},
then it suffices to focus on the asymptotics of $\pi_{T_u}(u)$.
We note that, independently of the value of $\varphi$,
$\frac{m(u)\Delta(u)}{u^{1/\beta}}\rw 0$ as $u\rw \IF$.

Let $D_{k}(u)=[k\frac{\Delta(u)}{u^{1/\beta}}S,(k+1)\frac{\Delta(u)}{u^{1/\beta}}S]$,
$F_{l}(u)=[\tau_u+l\frac{\Delta(u)}{u^{1/\beta}}S,\tau_u+(l+1)\frac{\Delta(u)}{u^{1/\beta}}S]$
and $I_{k,l}(u)=D_{k}(u)\times F_l(u)$ with $S>0$.
Moreover, let
$N_{S,u}=[\frac{u^{1/\beta}\ln m(u)}{m(u)\Delta(u)S}]$ and $m^{\pm \epsilon}_{k,l}(u)=m(u)\left(1+(b\pm \epsilon)\left((l-k)\frac{\Delta(u)}{u^{1/\beta}}S\right)^2\right)$.


{\bf Proof of case $\lim_{u\rw\IF}\frac{T_u}{\Delta(u)}=\IF$}.
\\
{\it \underline{Upper bound of $\pi_{T_u}(u)$}.}
Clearly, we have
\BQN\label{F1}
\pi_{T_u}(u)
&\leq&
\sum_{k=0}^{[\frac{T_u}{\Delta(u)S}]+1}\sum_{l=-N_{S,u}-1+k}^{N_{S,u}+2+k}\pl{\sup_{(t,s)\in I_{k,l}(u)}Z_u(s,t)>m(u)}\nonumber\\
&\leq& \sum_{k=0}^{[\frac{T_u}{\Delta(u)S}]+1}\sum_{l=-N_{S,u}-1+k}^{N_{S,u}+2+k}\pl{\sup_{(t,s)\in I_{k,l}(u)}\overline{Z_u}(s,t)>m^{-\epsilon}_{k,l}(u)},
\EQN
In order to apply Lemma \ref{PICKANDS}, we have to check conditions {\bf P1-P4},
for appropriately chosen $K_u, g_{k,l},\theta_{k,l}$.
Let
\[
X^{(u, k,l)}(t,s):=\overline{Z_u}\left(\tau_u+l\frac{\Delta(u)}{u^{1/\beta}}S+
\frac{\Delta(u)}{u^{1/\beta}}s, k\frac{\Delta(u)}{u^{1/\beta}}S+
\frac{\Delta(u)}{u^{1/\beta}}t\right)
\]
with $(s,t)\in[0,S]^2$ and $(k,l)\in K_u:=\{(k,l),
0\leq k\leq [\frac{T_u}{\Delta(u)S}]+1, -N_{S,u}-1+k\leq l \leq N_{S,u}+2+k\}$.
Then, let $g_{k,l}(u):=m^{-\epsilon}_{k,l}(u)$ and
\BQN\label{theta}
\theta_{k,l}(u,t,s,t_1,s_1)
:=\frac{\sigma^2(\Delta(u)|s-s_1|)+\sigma^2(\Delta(u)|t-t_1|)}{\sigma^2(\Delta(u))}
\frac{\sigma^2(\Delta(u))}{2\sigma^2(u^{1/\beta}\tau^*)}(m^{-\epsilon}_{k,l}(u))^2,
\EQN
for $(t,s), (t_1,s_1)\in[0,S]^2, (k,l)\in K_u.$

Assumption {\bf P1} holds straightforwardly. In order to show {\bf P2}, we observe that,
by definition of $\Delta(u)$,
\BQN\label{the1}
\lim_{u\rw\IF}\sup_{(k,l)\in K_u}\left|\frac{\sigma^2(\Delta(u))}{2\sigma^2(u^{1/\beta}\tau^*)}(m^{-\epsilon}_{k,l}(u))^2-1\right|=0
\EQN
and, by the UCT,
\BQNY
\lim_{u\rw\IF}\left|\frac{\sigma^2(\Delta(u)|s-s_1|)+\sigma^2(\Delta(u)|t-t_1|)}{\sigma^2(\Delta(u))}-|s-s_1|^{2{\alpha_0}}-|t-t_1|^{2{\alpha_0}}\right|=0, \ \
(t,s), (t_1,s_1)\in[0,S]^2.
\EQNY
Therefore using UCT again, we conclude that
\BQNY
&&\lim_{u\rw\IF}\sup_{(k,l)\in K_u}\left|\theta_{k,l}(u,t,s,t_1,s_1)-|s-s_1|^{2\alpha_0}-|t-t_1|^{2\alpha_0}\right|\\
&\leq& \lim_{u\rw\IF}\left|\frac{\sigma^2(\Delta(u)|s-s_1|)+\sigma^2(\Delta(u)|t-t_1|)}{\sigma^2(\Delta(u))}-|s-s_1|^{2\alpha_0}-|t-t_1|^{2\alpha_0}\right|\\
&&+\lim_{u\rw\IF}\sup_{(k,l)\in K_u}\left|\frac{\sigma^2(\Delta(u))}{2\sigma^2(u^{1/\beta}\tau^*)}(m^{-\epsilon}_{k,l}(u))^2-1\right|
\frac{\sigma^2(\Delta(u)|s-s_1|)+\sigma^2(\Delta(u)|t-t_1|)}{\sigma^2(\Delta(u))}\\
&\leq&\lim_{u\rw\IF}\left|\frac{\sigma^2(\Delta(u)|s-s_1|)+\sigma^2(\Delta(u)|t-t_1|)}{\sigma^2(\Delta(u))}-|s-s_1|^{2\alpha_0}-|t-t_1|^{2\alpha_0}\right|\\
&&+\lim_{u\rw\IF}\sup_{(k,l)\in K_u}\mathbb{Q}S^{2\alpha_0}\left|\frac{\sigma^2(\Delta(u))}{2\sigma^2(u^{1/\beta}\tau^*)}(m^{-\epsilon}_{k,l}(u))^2-1\right|\\
&\rw& 0, \ \ (t,s), (t_1,s_1)\in[0,S]^2,
\EQNY
which implies that {\bf P2} is satisfied.

In order to check {\bf P3}, we use that by UCT, with $\gamma_1$ and $g_1(t)$ defined
in (\ref{EG}),
\BQNY
&&\overline{\lim}_{u\rw\IF}\sup_{(k,l)\in K_u}\sup_{(t,s)\neq(t_1,s_1)\in[0,S]^2}\frac{\theta_{k,l}(u,t,s,t_1,s_1)}{|s-s_1|^{\gamma_1}+|t-t_1|^{\gamma_1}}\\
&\leq& \overline{\lim}_{u\rw\IF}\sup_{(k,l)\in K_u}\sup_{(t,s)\neq (t_1,s_1)\in[0,S]^2}2\frac{\sigma^2(\Delta(u)|s-s_1|)+\sigma^2(\Delta(u)|t-t_1|)}{\sigma^2(\Delta(u))(|s-s_1|^{\gamma_1}+|t-t_1|^{\gamma_1})}\\
&\leq & 2\overline{\lim}_{u\rw\IF}\sup_{(t,s)\neq (t_1,s_1)\in[0,S]^2}\frac{\sigma^2(\Delta(u)|s-s_1|)}{\sigma^2(\Delta(u))|s-s_1|^{\gamma_1}}+2\overline{\lim}_{u\rw\IF}
\sup_{(t,s)\neq (t_1,s_1)\in[0,S]^2}
\frac{\sigma^2(\Delta(u)|t-t_1|)}{\sigma^2(\Delta(u))|t-t_1|^{\gamma_1}}\\
&=&4\overline{\lim}_{u\rw\IF}\sup_{t\in[0,S]}\frac{g_1(\Delta(u)t)}{g_1(\Delta(u))}\leq 8S^{2\alpha_0-\gamma_1}<\IF.
\EQNY
Next we focus on {\bf P4}. First, in light of UCT and (\ref{the1})
\BQNY
&&\lim_{\epsilon\rw 0}\overline{\lim}_{u\rw\IF}\sup_{(k,l)\in K_u}\sup_{|(t,s)-(t_1,s_1)|<\epsilon, (t,s),(t_1,s_1)\in[0,S]^2}\left|\theta_{k,l}(u,t,s,0,0)-\theta_{k,l}(u,t_1,s_1,0,0)\right|\\
&\leq&2\lim_{\epsilon\rw 0}\overline{\lim}_{u\rw\IF}\sup_{s,t,s_1,t_1\in[0,\epsilon], }\left|\frac{\sigma^2(\Delta(u)s)+\sigma^2(\Delta(u)t)-\sigma^2(\Delta(u)s_1)-\sigma^2(\Delta(u)t_1)}{\sigma^2(\Delta(u))}\right|=0.
\EQNY
Second, it follows from Lemma \ref{L2} and and UCT that
\BQNY
&&\left|(m^{-\epsilon}_{k,l}(u))^2(1-r_u(s_l(u)+\frac{\Delta(u)}{u^{1/\beta}}s,t_k(u)+\frac{\Delta(u)}{u^{1/\beta}}t, s_l(u),t_k(u)))-\theta_{k,l}(u,t,s,0,0)\right|\\
&\leq &\left|(m^{-\epsilon}_{k,l}(u))^2\frac{1-r_u(s_l(u)+\frac{\Delta(u)}{u^{1/\beta}}s,t_k(u)+\frac{\Delta(u)}{u^{1/\beta}}t, s_l(u),t_k(u))}{\theta_{k,l}(u,t,s,0,0)}-1\right|\theta_{k,l}(u,t,s,0,0)\\
&\leq& \mathbb{Q}_0S^{2\alpha_0}\left|(m^{-\epsilon}_{k,l}(u))^2\frac{1-r_u(s_l(u)+\frac{\Delta(u)}{u^{1/\beta}}s,t_k(u)+\frac{\Delta(u)}{u^{1/\beta}}t, s_l(u),t_k(u))}{\theta_{k,l}(u,t,s,0,0)}-1\right|\RW 0, \ \ u\rw\IF,
\EQNY
with respect to $(k,l)\in K_u$ and $(t,s)\in[0,S]^2$, where $s_l(u)=\tau_u+l\frac{\Delta(u)}{u^{1/\beta}}S$ and $t_k(u)=k\frac{\Delta(u)}{u^{1/\beta}}S$.

The above leads to, for $\|(t,s)-(t_1,s_1)\|<\epsilon$,
\BQNY
&&(m^{-\epsilon}_{k,l}(u))^2\E{\left(X^{(u,k,l)}(s,t)-X^{(u,k,l)}(s_1,t_1)\right)X^{(u,k,l)}(0,0)}\\
&\leq& \left|(m^{-\epsilon}_{k,l}(u))^2(1-r_u(s_l(u)+\frac{\Delta(u)}{u^{1/\beta}}s,t_k(u)+\frac{\Delta(u)}{u^{1/\beta}}t, s_l(u),t_k(u)))-\theta_{k,l}(u,t,s,0,0)\right|\\
&&+\left|(m^{-\epsilon}_{k,l}(u))^2(1-r_u(s_l(u)+\frac{\Delta(u)}{u^{1/\beta}}s_1,t_k(u)+\frac{\Delta(u)}{u^{1/\beta}}t_1, s_l(u),t_k(u)))-\theta_{k,l}(u,t_1,s_1,0,0)\right|\\
&&+\left|\theta_{k,l}(u,t,s,0,0)-\theta_{k,l}(u,t_1,s_1,0,0)\right|\RW 0,\ \ u\rw\IF, \epsilon\rw 0.
\EQNY
with respect to $(k,l)\in K_u, (t,s),(t_1,s_1)\in[0,S]^2$,
which confirms that {\bf P4} is fulfilled.

Thus in view of Lemma \ref{PICKANDS}
\BQN\label{PI_1}
\frac{\pl{\sup_{(t,s)\in I_{k,l}(u)}\overline{Z_u}(s,t)>m^{-\epsilon}_{k,l}(u)}}{\Psi(m^{-\epsilon}_{k,l}(u))}\to \mathcal{H}_V([0,S]^2)= \left(\mathcal{H}_{B_{\alpha_0}}[0,S]\right)^2,\ \ u\rw\IF,
\EQN
where $V(t,s)=B^{(1)}_{{\alpha_0}}(t)+B^{(2)}_{{\alpha_0}}(s)$,
with $B^{(1)}_{{\alpha_0}}$ and $B^{(2)}_{{\alpha_0}}$ being
two independent fBms with index ${\alpha_0}$. Then, continuing (\ref{F1}), in view of (\ref{PB}) we have
\BQNY
\pi_{T_u}(u)
&\leq&\sum_{k=0}^{[\frac{T_u}{\Delta(u)S}]+1}\sum_{l=-N_{S,u}-1+k}^{N_{S,u}+2+k}(\mathcal{H}_{B_{\alpha_0}}[0,S])^2\Psi(m_{k,l}^{-\epsilon}(u))(1+o(1)) \\
&\leq&\sum_{k=0}^{[\frac{T_u}{\Delta(u)S}]+1}(\mathcal{H}_{B_{\alpha_0}}[0,S])^2\Psi(m(u))\sum_{l=-N_{S,u}-2}^{N_{S,u}+3}e^{-(b-\epsilon)\left(lm(u)
\frac{\Delta(u)}{u^{1/\beta}}S\right)^2}(1+o(1))\nonumber\\
&\leq& \left(\frac{\mathcal{H}_{B_{{\alpha_0}}}[0,S]}{S}\right)^2(b-\epsilon)^{-1/2}[\frac{T_u}{\Delta(u)}]\frac{u^{1/\beta}}{m(u)\Delta(u)}\Psi(m(u))
\int_{-\IF}^{\IF}e^{-x^2}dx(1+o(1))\nonumber\\
&\sim& (\mathcal{H}_{B_{{\alpha_0}}})^2(b-\epsilon)^{-1/2}\sqrt{\pi}\frac{u^{1/\beta}T_u}{(\Delta(u))^2m(u)}\Psi(m(u)), \ \ \mbox{as}\ \ u\rw\IF,
\EQNY
with $b=\frac{B}{2A}$ (see Lemma \ref{L1}).
Hence, letting $\epsilon\rw 0$, we obtain the upper bound for $\pi_{T_u}(u)$.
\\
{\it \underline{Lower bound of $\pi_{T_u}(u)$}}.
Set
\\
$\Gamma_{\delta,1}=\{(k,l,,k_1,l_1):
0\leq k\leq k_1\leq [\frac{T_u}{\Delta(u)S}],
|k_1-k|\leq \frac{\delta u^{1/\beta}}{\Delta(u)S},
-N_{S,u}+k\leq l\leq l_1\leq N_{S,u}+k, I_{k,l}(u)\cap
I_{k_1,l_1}(u)=\emptyset\}$,
$\Gamma_{\delta,2}=\{(k,l,,k_1,l_1): 0\leq k\leq k_1\leq [\frac{T_u}{\Delta(u)S}], |
k_1-k|\leq \frac{\delta u^{1/\beta}}{\Delta(u)S},
-N_{S,u}+k\leq l\leq l_1\leq N_{S,u}+k, I_{k,l}(u)\cap I_{k_1,l_1}(u)\neq \emptyset\}$,
$\Gamma_{\delta,3}=\{(k,l,,k_1,l_1): 0\leq k\leq k_1\leq [\frac{T_u}{\Delta(u)S}],
\frac{\delta u^{1/\beta}}{\Delta(u)S}<|k_1-k|\leq\frac{u^{1/\beta}}
{\Delta(u)S}e^{\frac{1-a_\delta}{4+4a_\delta}m^2(u)}, -N_{S,u}+k\leq l\leq l_1\leq N_{S,u}+k \}$,
$\Gamma_{\delta,4}=\{(k,l,,k_1,l_1): 0\leq k\leq k_1\leq [\frac{T_u}{\Delta(u)S}], |k_1-k|> \frac{u^{1/\beta}}{\Delta(u)S}e^{\frac{1-a_\delta}{4+4a_\delta}m^2(u)}, -N_{S,u}+k\leq l\leq l_1\leq N_{S,u}+k \}$.
We have
 \BQNY
 \pi_{T_u}(u)\geq\sum_{k=0}^{[\frac{T_u}{\Delta(u)S}]}\sum_{l=-N_{S,u}+k}^{N_{S,u}+k}\pl{\sup_{(t,s)\in I_{k,l}(u)}Z_u(s,t)>m(u)}
 -\left(\Sigma_1(u)+\Sigma_2(u)+\Sigma_3(u)+\Sigma_4(u)\right),
 \EQNY
 where
 \BQNY
  \Sigma_i(u)&=&\sum_{(k,l,k_1,l_1)\in\Gamma_{\delta,i}}\pl{\sup_{(t,s)\in I_{k,l}(u)}Z_u(s,t)>m(u)\sup_{(t_1,s_1)\in I_{k_1,l_1}(u)}Z_u(s_1,t_1)>m(u)}, \ \ i=1,2,3,4.
 \EQNY
 The same lines of reasoning, as presented in the
 proof of the upper bound of $\pi_{T_u}(u)$,
give the lower bound
for
$\sum_{k=0}^{[\frac{T_u}{\Delta(u)S}]}\sum_{l=-N_{S,u}+k}^{N_{S,u}+k}\pl{\sup_{(t,s)\in I_{k,l}(u)}Z_u(s,t)>m(u)}$,
which asymptotically agrees with the upper bound.
Thus the remaining task is to prove that
$\Sigma_i(u)$, $i=1,2,3,4$ are asymptotically negligible.

{\it \underline{Upper bound of $\Sigma_1(u)$}}.
In light of Lemma \ref{L2}, there exists a positive constant
$\delta>0$ such that for $u$ large enough,
all $(t,s,t_1,s_1)\in I_{k,l}(u)\times I_{k_1,l_1}(u)$
with $(k,l,k_1,l_1)\in \Gamma_{\delta,1}$,
 \BQNY
 1/2<\frac{1-r_u(s,t,s_1,t_1)}{\frac{\sigma^2(u^{1/\beta}|s-s_1|)+\sigma^2(u^{1/\beta}|t-t_1|)}{2\sigma^2(u^{1/\beta}\tau^*)}}<2.
 \EQNY
Moreover, by UCT, we have
 \BQNY
 2\leq \E{\overline{Z_u}(s,t)+\overline{Z_u}(s_1,t_1)}^2&=&4-2(1-r_u(s,t,s_1,t_1))\\
 &\leq& 4-\frac{\sigma^2(u^{1/\beta}|s-s_1|)+\sigma^2(u^{1/\beta}|t-t_1|)}{2\sigma^2(u^{1/\beta}\tau^*)}\\
 &\leq& 4-\mathbb{Q}_2\frac{|l_1-l|^{\gamma_1}S^{\gamma_1}+|k_1-k|^{\gamma_1}S^{\gamma_1}}{m^2(u)},
 \EQNY
 where $0<\gamma_1<\min(2{\alpha_\IF},\gamma)$.

Thus
 \BQNY
 \Sigma_1(u)&\leq& \sum_{(k,l,k_1,l_1)\in\Gamma_{\delta,1}}\pl{\sup_{(t,s)\in I_{k,l}(u)}\overline{Z_u}(s,t)>m_{k,l}^{-\epsilon}(u),\sup_{(t_1,s_1)\in I_{k_1,l_1}(u)}\overline{Z_u}(s_1,t_1)>m_{k_1,l_1}^{-\epsilon}(u)}\\
 &\leq& \sum_{(k,l,k_1,l_1)\in\Gamma_{\delta,1}}\pl{\sup_{(t,s,t_1,s_1)\in I_{k,l}(u)\times I_{k_1,l_1}(u)}\overline{Z_u}(s,t)+\overline{Z_u}(s_1,t_1)>2\hat{m}_{k,l,k_1,l_1}^{-\epsilon}(u)}\\
 &\leq & \sum_{(k,l,k_1,l_1)\in\Gamma_{\delta,1}}\pl{\sup_{(t,s,t_1,s_1)\in I_{k,l}(u)\times I_{k_1,l_1}(u)}\overline{\overline{Z_u}(s,t)+\overline{Z_u}(s_1,t_1)}>\frac{2\hat{m}_{k,l,k_1,l_1}^{-\epsilon}(u)}
 {\sqrt{4-\mathbb{Q}_2\frac{|l_1-l|^{\gamma_1}
 S^{\gamma_1}+|k_1-k|^{\gamma_1}S^{\gamma_1}}{m^2(u)}}}}
 \EQNY
with $\hat{m}_{k,l,k_1,l_1}^{-\epsilon}(u)=\min(m_{k,l}^{-\epsilon}(u),m_{k_1,l_1}^{-\epsilon}(u))$.\\
In order to bound the above sum, we introduce
\[r_u(t,s,t_1,s_1,t',s',t_1^{'},s_1^{'})
:=\E{(\overline{\overline{Z_u}(s,t)+\overline{Z_u}(s_1,t_1)})
(\overline{\overline{Z_u}(s',t')+\overline{Z_u}(s_1^{'},t_1^{'})})}\]
and observe that for
$(t,s,t_1,s_1)$, $(t',s',t_1^{'},s_1^{'})$ $\in I_{k,l}(u)\times I_{k_1,l_1}(u)$,
\BQNY
1-r_u(t,s,t_1,s_1,t',s',t_1^{'},s_1^{'})&\leq& \frac{\E{\overline{Z_u}(s,t)+\overline{Z_u}(s_1,t_1)-\overline{Z_u}(s',t')-\overline{Z_u}(s_1^{'},t_1^{'})}^2}
{2\sqrt{\E{\overline{Z_u}(s,t)+\overline{Z_u}(s_1,t_1)}^2}\sqrt{\E{\overline{Z_u}(s',t')+\overline{Z_u}(s_1^{'},t_1^{'})}^2}}\\
&\leq&\frac{\E{\overline{Z_u}(s,t)-\overline{Z_u}(s',t')}^2+\E{\overline{Z_u}(s_1,t_1)-\overline{Z_u}(s_1^{'},t_1^{'})}^2}{2}\\
&\leq& 1-r_u(s,t,s',t')+1-r_u(s_1,t_1,s_1^{'},t_1^{'})\\
&\leq&\frac{\sigma^2(u^{1/\beta}|s-s^{'}|)+\sigma^2(u^{1/\beta}|t-t^{'}|)}{\sigma^2(u^{1/\beta}\tau^*)}+
\frac{\sigma^2(u^{1/\beta}|s_1-s^{'}_1|)+\sigma^2(u^{1/\beta}|t_1-t^{'}_1|)}{\sigma^2(u^{1/\beta}\tau^*)}\\
&\leq&\mathbb{Q}_3S^2\frac{\left(\frac{u^{1/\beta}}{\Delta(u)}\right)^{\gamma_2}
\left(|s-s'|^{\gamma_2}+|t-t'|^{\gamma_2}+|s_1-s_1^{'}|^{\gamma_2}+|t_1-t_1^{'}|^{\gamma_2}\right)}{m^2(u)},
\EQNY
with $0<\gamma_2<\min(2{\alpha_\IF},\gamma)$ and $S\geq 1$.

Next we define a centered homogenous Gaussian field
$\{X_u^*(s,t,s_1,t_1), (s,t,s_1,t_1)\in \mathbb{R}^4\}$
so that
$X_u^*(s,t,s_1,t_1):=(X^{1}_u(s)+X^{2}_u(t)+X^{3}_u(s_1)+X^{4}_u(t_1))/2$
with $X_u^{i}(s), 1\leq i\leq 4,$ being i.i.d. centered stationary Gaussian processes
with covariance function
\BQNY
r_u(s,s')=\exp\left(-8\mathbb{Q}_3S^2\left(\frac{u^{1/\beta}}{\Delta(u)}\right)^{\gamma_2}
\frac{1}{m^2(u)}|s-s'|^{\gamma_2}\right).
\EQNY
Let $r_u^{*}(s,t,s_1,t_1,s',t',s_1^{'},t_1^{'})$ be the covariance function of $X_u^*(s,t,s_1,t_1)$.
It is straightforward to check that
for $(t,s,t_1,s_1)$, $(t',s',t_1^{'},s_1^{'})$ $\in I_{k,l}(u)\times I_{k_1,l_1}(u)$,
\BQNY
r_u(s,t,s_1,t_1,s',t',s_1^{'},t_1^{'})
\ge
r_u^{*}(s,t,s_1,t_1,s',t',s_1^{'},t_1^{'}).
\EQNY
In light of Slepian's inequality ( see, e.g., \cite{AdlerTaylor} or \cite{Pit96}), we have
\BQN\label{F2}
\Sigma_1(u)&\leq& \sum_{(k,l,k_1,l_1)\in\Gamma_{\delta,1}}\pl{\sup_{(t,s,t_1,s_1)\in I_{k,l}(u)\times I_{k_1,l_1}(u)}X^*_u(s,t,s_1,t_1)>\frac{2\hat{m}_{k,l,k_1,l_1}^{-\epsilon}(u)}{\sqrt{4-\mathbb{Q}_2\frac{|l_1-l|^{\gamma_1}
 S^{\gamma_1}+|k_1-k|^{\gamma_1}S^{\gamma_1}}{m^2(u)}}}}\nonumber\\
 &\leq&\sum_{(k,l,k_1,l_1)\in\Gamma_{\delta,1}}2(\mathcal{H}_{B_{{\gamma_2}}}[0,S_1])^4\Psi\left(\frac{2\hat{m}_{k,l,k_1,l_1}^{-\epsilon}(u)}
 {\sqrt{4-\mathbb{Q}_2\frac{|l_1-l|^{\gamma_1}
 S^{\gamma_1}+|k_1-k|^{\gamma_1}S^{\gamma_1}}{m^2(u)}}}\right)\nonumber\\
 &\leq&\sum_{(k,l,k_1,l_1)\in\Gamma_{\delta,1}}4(\mathcal{H}_{B_{\gamma_2}}[0,S_1])^4\Psi(\hat{m}_{k,l,k_1,l_1}^{-\epsilon}(u))e^{-\mathbb{Q}_4\left(|l_1-l|^{\gamma_1}
 S^{\gamma_1}+|k_1-k|^{\gamma_1}S^{\gamma_1}\right)}\nonumber\\
 &\leq&\sum_{k=0}^{[\frac{T_u}{\Delta(u)S}]}\sum_{l=-N_{S,u}+k}^{N_{S,u}+k}4\left(\frac{\mathcal{H}_{B_{{\gamma_2}}}[0,S_1]}{S_1}\right)^4
 \Psi(\hat{m}_{k,l}^{-\epsilon}(u))S_1^{-4}\sum_{i\geq 0,j\geq 0,i+j\geq 1}e^{-\mathbb{Q}_4\left(i^{\gamma_1}S^{\gamma_1}+j^{\gamma_1}S^{\gamma_1}\right)}\nonumber\\
 &\leq&8\left(\frac{\mathcal{H}_{B_{\gamma_2}}[0,S_1]}{S_1}\right)^4(b-\epsilon)^{-1/2}\sqrt{\pi}\frac{u^{1/\beta}T_u}{(\Delta(u))^2m(u)}
 \Psi(m(u))S_1^{-4}e^{-\mathbb{Q}_5S^{\gamma_1}}
\EQN
with $S_1=(2\mathbb{Q}_3)^{2/\gamma_2}S^{1+4/\gamma_2}$.
Letting $S\rw \IF$, we get that $\Sigma_1(u)=o(\pi_{T_u}(u))$
as $u\to\infty$.

{\it \underline{Upper bound of $\Sigma_2(u)$}}.
We have
\begin{eqnarray*}
\Sigma_2(u)\le
\sum_{(k,l,k_1,l_1)\in\Gamma_{\delta,2}}p^{(1)}_{k,l, k_1,l_1}(u)+
\sum_{(k,l,k_1,l_1)\in\Gamma_{\delta,2}}p^{(2)}_{k,l, k_1,l_1}(u),
\end{eqnarray*}
with
\[
p^{(1)}_{k,l,k_1,l_1}(u)
:=
\pl{\sup_{(t,s)\in I_{k,l}(u)}Z_u(s,t)>m(u) \sup_{(t_1,s_1)\in I_{k_1,l_1}^{1}(u)}Z_u(s_1,t_1)>m(u)}
\]
\[
p^{(2)}_{k,l,k_1,l_1}(u)
:=
\pl{\sup_{(t,s)\in I_{k,l}(u)}Z_u(s,t)>m(u)\sup_{(t_1,s_1)\in I_{k_1,l_1}^{2}(u)}Z_u(s_1,t_1)>m(u)},
\]
where, without loss of generality we assume that
$k+1=k_1$ and, $l=l_1$ or $l\pm1=l_1$
and
\[I_{k_1,l_1}^{1}(u)=[(k+1)\frac{\Delta(u)}{u^{1/\beta}}S,(k+1)\frac{\Delta(u)}{u^{1/\beta}}S+\frac{\Delta(u)}{u^{1/\beta}}\sqrt{S}]\times F_{l_1}(u),\]
\[I_{k_1,l_1}^{2}(u)=[(k+1)\frac{\Delta(u)}{u^{1/\beta}}S+\frac{\Delta(u)}{u^{1/\beta}}\sqrt{S},(k+2)\frac{\Delta(u)}{u^{1/\beta}}S,]\times F_{l_1}(u).\]
Following the same argument as given in the proof of $\Sigma_1(u)$,
we get
\BQN\label{E1}
\sum_{(k,l,k_1,l_1)\in\Gamma_{\delta,2}}p^{(1)}_{k,l, k_1,l_1}(u)&\leq& \sum_{k_1=0}^{[\frac{T_u}{\Delta(u)S}]+1}\sum_{l_1=-N_{S,u}-1+k}^{N_{S,u}+2+k}\pl{\sup_{(t_1,s_1)\in I_{k_1,l_1}^{1}(u)}Z_u(s_1,t_1)>m(u)}\nonumber
 \\ &\leq&\sum_{k=0}^{[\frac{T_u}{\Delta(u)S}]+1}\sum_{l=-N_{S,u}-1+k}^{N_{S,u}+2+k}\mathcal{H}_{B_{\alpha_0}}[0,\sqrt{S}]\mathcal{H}_{B_{\alpha_0}}[0,S]
\Psi(m_{k,l}^{-\epsilon}(u))(1+o(1))\nonumber\\
&\leq&\frac{2}{\sqrt{S}}\frac{\mathcal{H}_{B_{\alpha_0}}[0,\sqrt{S}]}{\sqrt{S}}\frac{\mathcal{H}_{B_{\alpha_0}}[0,S]}{S}(b-\epsilon)^{-1/2}
\sqrt{\pi}\frac{u^{1/\beta}T_u}{(\Delta(u))^2m(u)}\Psi(m(u)).
\EQN
and, with the same $S_1$ as above,
\BQN\label{E2}
\sum_{(k,l,k_1,l_1)\in\Gamma_{\delta,2}}p^{(2)}_{k,l, k_1,l_1}(u)
&&\leq \sum_{(k,l,k_1,l_1)\in\Gamma_{\delta,2}}\pl{\sup_{(t,s,t_1,s_1)
\in I_{k,l}(u)\times I_{k_1,l_1}(u)}\overline{\overline{Z_u}(s,t)+\overline{Z_u}(s_1,t_1)}>\frac{2\hat{m}_{k,l,k_1,l_1}^{-\epsilon}(u)}{\sqrt{4-\mathbb{Q}_3\frac{
 S^{\gamma_1/2}}{m^2(u)}}}}\nonumber\\
 &&\leq\sum_{(k,l,k_1,l_1)\in\Gamma_{\delta,2}}2(\mathcal{H}_{B_{{\alpha_0}}}[0,S_1])^4\Psi\left(\frac{2\hat{m}_{k,l,k_1,l_1}^{-\epsilon}(u)}
 {\sqrt{4-\mathbb{Q}_3\frac{
 S^{\gamma_1/2}}{m^2(u)}}}\right)\nonumber\\
 &&\leq\sum_{k=0}^{[\frac{T_u}{\Delta(u)S}]}\sum_{l=-N_{S,u}+k}^{N_{S,u}+k}4\left(\frac{\mathcal{H}_{B_{{\alpha_0}}}[0,S_1]}{S_1}\right)^4
 \Psi(\hat{m}_{k,l}^{-\epsilon}(u))S_1^{-4}
 e^{-\mathbb{Q}_3S^{\gamma_1/2}}\nonumber\\
 &&\leq 8\left(\frac{\mathcal{H}_{B_{{\alpha_0}}}[0,S_1]}{S_1}\right)^4(b-\epsilon)^{-1/2}\sqrt{\pi}\frac{u^{1/\beta}T_u}{(\Delta(u))^2m(u)}\Psi(m(u))S_1^{-4}
 e^{-\mathbb{Q}_3S^{\gamma_1/2}}.
\EQN
Combination of (\ref{E1}) with (\ref{E2}) implies that
$\Sigma_2(u)=o(\pi_{T_u}(u))$
as $u\to\infty$.

{\it \underline{Upper bound of $\Sigma_3(u)$}}.
The idea of this part of the proof is to apply Borel inequality.
For that, without loss of generality, we fix $S=1$.
We observe that (similarly as for $\Sigma_1(u)$),
for
$(t,s,t_1,s_1), (t',s',t^{'}_1,s^{'}_1)\in I_{k,l}(u)\times I_{k_1,l_1}(u)$
with $(k,l,k_1,l_1)\in \Gamma_{\delta,3}$,
\BQNY
\E{\overline{Z_u}(s,t)+\overline{Z_u}(s_1,t_1)}^2=2+2r_u(s,t,s_1,t_1)\leq 2+2a_\delta<4
\EQNY
and
\BQNY
&&\E{\overline{Z_u}(s,t)+\overline{Z_u}(s_1,t_1)-\overline{Z_u}(s',t')-\overline{Z_u}(s^{'}_1,t^{'}_1)}^2\\
&\leq& 4(1-r_u(s,t,s',t'))+4(1-r_u(s_1,s_1^{'},t_1,t_1^{'}))\\
&\leq&\mathbb{Q}_6\frac{\left(\frac{u^{1/\beta}}{\Delta(u)}\right)^{\gamma_2}
\left(|s-s'|^{\gamma_2}+|t-t'|^{\gamma_2}+|s_1-s_1^{'}|^{\gamma_2}+|t_1-t_1^{'}|^{\gamma_2}\right)}{m^2(u)}.
\EQNY
Thus, by Fernique inequality,
\BQNY
\pl{\sup_{(t,s,t_1,s_1)\in I_{k,l}(u)\times I_{k_1,l_1}(u)}\overline{Z_u}(s,t)+\overline{Z_u}(s_1,t_1)>x}\leq \frac{1}{2} e^{-\frac{x^2}{8}},
\EQNY
for any $(k,l,k_1,l_1)\in\Gamma_{\delta,3}$, any $x>0$ and $u$ large enough.
This implies that there exists a common positive constant $a$ such that
for any $(k,l,k_1,l_1)\in\Gamma_{\delta,3}$ and $u$ large enough
\BQNY
\pl{\sup_{(t,s,t_1,s_1)\in I_{k,l}(u)\times I_{k_1,l_1}(u)}\overline{Z_u}(s,t)+\overline{Z_u}(s_1,t_1)>a}\leq 1/2.
\EQNY
The above implies that we can apply Borel inequality to the sum below uniformly
\BQNY
\Sigma_3(u)&\leq& \sum_{(k,l,k_1,l_1)\in\Gamma_{\delta,3}}\pl{\sup_{(t,s)\in I_{k,l}(u)}\overline{Z_u}(s,t)>m(u),\sup_{(t_1,s_1)\in I_{k_1,l_1}(u)}\overline{Z_u}(s_1,t_1)>m(u)},\\
&\leq& \sum_{(k,l,k_1,l_1)\in\Gamma_{\delta,3}}\pl{\sup_{(t,s,t_1,s_1)\in I_{k,l}(u)\times I_{k_1,l_1}(u)}\overline{Z_u}(s,t)+\overline{Z_u}(s_1,t_1)>2m(u)}\\
&\leq & \mathbb{Q}_7\left(\frac{u^{1/(2\beta)}}{\Delta(u)}\frac{u^{1/\beta}\ln m(u)}{m(u)\Delta(u)}\right)^2T_ue^{\frac{1-a_\delta}{4+4a_\delta}m^2(u)}\Psi\left(\frac{2m(u)-a}{\sqrt{2+2a_\delta}}\right)\\
 &\leq& \mathbb{Q}_8\left(\frac{u^{3/(2\beta)}\ln m(u)}{m(u)\Delta^2(u)}\right)^2T_u\Psi\left(m(u)\right)e^{-m^2(u)\left(\frac{1-a_\delta}{4+4a_\delta}\right)+\frac{am(u)}{1+a_\delta}}\\
 &=&o\left(\frac{u^{1/\beta}T_u}{(\Delta(u))^2m(u)}\Psi(m(u))\right), \ \ u\rw\IF.
 \EQNY
This implies that $\Sigma_3(u)=o(\pi_{T_u}(u))$
as $u\to\infty$.

{\it \underline{Upper bound of $\Sigma_4(u)$}}.
Let $0<\epsilon<\frac{1-2\beta_1}{1+2\beta_1}$ be given.
Then, for $u$ large enough, $r_u(s,t,s_1,t_1)<\epsilon$ holds
for $|t-t_1|>e^{\frac{1-a_\delta}{8+8a_\delta}m^2(u)}$ and $s-t,s_1-t_1\in E(u)$.
Thus similarly as for $\Sigma_3(u)$, we have
 \BQNY
  \Sigma_4(u)&\leq& \mathbb{Q}_9\left(\frac{T_u}{\Delta(u)}\frac{u^{1/\beta}\ln m(u)}{m(u)\Delta(u)}\right)^2\Psi\left(\frac{2m(u)-a}{\sqrt{2+2\epsilon}}\right)\\
  &\leq& \mathbb{Q}_{10}\left(\frac{u^{1/\beta}\ln m(u)}{m(u)\Delta^2(u)}\right)^2T_u\Psi\left(m(u)\right)e^{-\left(\frac{1-\epsilon}{2(1+\epsilon)}-\beta_1\right)m^2(u)+\frac{am(u)}{1+\epsilon}}\\
   &=&o\left(\frac{u^{1/\beta}T_u}{(\Delta(u))^2m(u)}\Psi(m(u))\right), \ \ u\rw\IF.
 \EQNY
Hence $\Sigma_3(u)=o(\pi_{T_u}(u))$
as $u\to\infty$.

Note that if for some $T_u$,
$\Gamma_{\delta,3}$ or $\Gamma_{\delta,4}$
are empty then the above inequalities are still valid.
This completes the proof of ii).

{\bf Proof of case $\lim_{u\rw\IF}\frac{T_u}{\Delta(u)}=\rho\in(0,\IF)$}.
The proof of this case is similar to the proof of the previous case.
Thus we focus on the tiny details that differ from the arguments used in the previous case.

For $I_{0,l}^{(\pm\epsilon)}(u)=[0, (\rho\pm\epsilon)\frac{\Delta(u)}{u^{1/\beta}}]\times[\tau_u+l\frac{\Delta(u)}{u^{1/\beta}}S,\tau_u+(l+1)\frac{\Delta(u)}{u^{1/\beta}}S]$
we have
 \BQN\label{U}
\pi_{T_u}(u)
&\leq&
\sum_{l=-N_{S,u}-1}^{N_{S,u}+1}\pl{\sup_{(t,s)\in I_{0,l}^{+\epsilon}(u)}Z_u(s,t)>m(u)}\nonumber\\
&\leq& \sum_{l=-N_{S,u}-1}^{N_{S,u}+1}\pl{\sup_{(t,s)\in I_{0,l}^{+\epsilon}(u)}\overline{Z_u}(s,t)>m^{-\epsilon}_{0,l}(u)}\nonumber\\
&\leq& \sum_{l=-N_{S,u}-1}^{N_{S,u}+1}\mathcal{H}_{B_{\alpha_0}}[0,\rho+\epsilon]\mathcal{H}_{B_{\alpha_0}}[0,S]\Psi(m_{0,l}^{-\epsilon}(u))(1+o(1))\nonumber\\
&\sim& \mathcal{H}_{B_{{\alpha_0}}}[0,\rho+\epsilon]\mathcal{H}_{B_{{\alpha_0}}}(b-\epsilon)^{-1/2}\sqrt{\pi}\frac{u^{1/\beta}}{\Delta(u)m(u)}\Psi(m(u)),\label{b.new1}
\EQN
with $N_{S,u}$ and $m^{-\epsilon}_{k,l}(u)$ defined
below (\ref{F1}).
Similarly, we have
\BQN
 \pi_{T_u}(u)
 &\geq&
 \sum_{l=-N_{S,u}}^{N_{S,u}}\pl{\sup_{(t,s)\in I_{0,l}^{-\epsilon}(u)}Z_u(s,t)>m(u)}
 -\sum_{i=1}^2\Sigma_i'(u)\nonumber\\
&\ge& \mathcal{H}_{B_{{\alpha_0}}}[0,\rho-\epsilon]\mathcal{H}_{B_{{\alpha_0}}}(b-\epsilon)^{-1/2}\sqrt{\pi}\frac{u^{1/\beta}}{\Delta(u)m(u)}\Psi(m(u))(1+o(1))
-\sum_{i=1}^2\Sigma_i'(u),\label{L}
 \EQN
 where
 \BQNY
  \Sigma_i'(u)&=&\sum_{(l,l_1)\in\Gamma_{i}'}\pl{\sup_{(t,s)\in I_{0,l}^{-\epsilon}(u)}Z_u(s,t)>m(u)\sup_{(t_1,s_1)\in I_{0,l_1}^{-\epsilon}(u)}Z_u(s_1,t_1)>m(u)}, \ \ i=1,2,
 \EQNY
 with $\Gamma_{1}'=\{(l,l_1), -N_{S,u}\leq l< l_1+1\leq N_{S,u} \}$ and $\Gamma_{2}'=\{(l,l_1), -N_{S,u}\leq l= l_1+1\leq N_{S,u} \}$.\\
 Following the same lines of argument
 as in (\ref{F2}) (see also (\ref{E1}) or (\ref{E2})),
 we get that $\sum_{i=1}^2\Sigma_i'(u)$ is negligible
 compared with the first term in (\ref{L}).
 Hence, comparing (\ref{U}) with (\ref{L}) and letting $\epsilon\rw 0$,
 we obtain that for $\frac{T_u}{\Delta(u)}\rw \rho\in (0,\infty)$,
 \BQN
 \Psi_{T_u}(u)\sim\mathcal{H}_{B_{{\alpha_0}}}[0,\rho]\mathcal{H}_{B_{{\alpha_0}}}b^{-1/2}\sqrt{\pi}\frac{u^{1/\beta}}{\Delta(u)m(u)}\Psi(m(u))
 \label{fin.1}.
 \EQN
Finally let us suppose that
$\frac{T_u}{\Delta(u)}\rw \rho=0$.
Clearly, for any $\epsilon>0$,
 \BQNY
 \Psi_0(u)\leq \Psi_{T_u}(u)\leq \Psi_{\Delta(u)\epsilon}(u).
 \EQNY
Hence, by (\ref{fin.1}),
$ \Psi_{\Delta(u)\epsilon}(u)\sim \mathcal{H}_{B_{{\alpha_0}}}[0,\epsilon]\mathcal{H}_{B_{{\alpha_0}}}b^{-1/2}\sqrt{\pi}\frac{u^{1/\beta}}{\Delta(u)m(u)}\Psi(m(u)).
$
Moreover, following \cite{DI2005}, we have
$ \Psi_0(u)\sim\mathcal{H}_{B_{{\alpha_0}}}b^{-1/2}\sqrt{\pi}\frac{u^{1/\beta}}{\Delta(u)m(u)}\Psi(m(u)).
$
Thus, using that
$\lim_{\epsilon\to 0} \mathcal{H}_{B_{{\alpha_0}}}[0,\epsilon]=1$,
we arrive at
 \BQNY
  \Psi_{T_u}(u)\sim\mathcal{H}_{B_{{\alpha_0}}}b^{-1/2}\sqrt{\pi}\frac{u^{1/\beta}}{\Delta(u)m(u)}\Psi(m(u)),
 \EQNY
which completes the proof.
\QED\\

\subsection{Proof of Theorem \ref{main4}}
In view of the proof of Theorem \ref{main3},
using the same notation for
$m^{\pm\epsilon}_{k,l}(u)$, $g_{k,l}(u)$ and $K_u$,
conditions
{\bf P1--P4} hold with
\BQN\label{theta1}
\theta_{k,l}(u,t,s,t_1,s_1)
:=\left(\sigma_1^2(|s-s_1|)+\sigma_1^2(|t-t_1|)\right)\frac{2\varphi^2(\tau^*)^{4{\alpha_\IF}}}{(1+c(\tau^*)^\beta)^2}
\frac{(m^{-\epsilon}_{k,l}(u))^2}{2\sigma^2(u^{1/\beta}\tau^*)},
\EQN
where $\sigma_1(t)=\frac{1+c(\tau^*)^\beta}{\sqrt{2}\varphi(\tau^*)^{2\alpha_\IF}}\sigma(t)$.
Thus, following Lemma \ref{PICKANDS},
\BQNY\label{PII}
\frac{\pl{\sup_{(t,s)\in I_{k,l}(u)}\overline{Z_u}(s,t)>m^{-\epsilon}_{k,l}(u)}}
{\Psi(m^{-\epsilon}_{k,l}(u))}\to \mathcal{H}_V([0,S]^2)= \left(\mathcal{H}_{\frac{1+c(\tau^*)^\beta}{\sqrt{2}\varphi(\tau^*)^{2\alpha_\IF}}X}[0,S]\right)^2,\ \ u\rw\IF,
\EQNY
where $V(t,s):=
\frac{1+c(\tau^*)^\beta}{\sqrt{2}\varphi(\tau^*)^{2\alpha_\IF}}
\left(X^{(1)}(t)+X^{(2)}(s)\right)$
with $X^{(1)},X^{(2)}$ being independent copies of $X$

The rest of the proof goes line-by-line the same as the proof of Theorem \ref{main3}.
\QED

\subsection{ Proof of Theorem \ref{main5}}
Similarly to the proof of Theorem \ref{main4},
{\bf P1-P4} hold with
\BQNY
\theta_{k,l}(u,t,s,t_1,s_1)
:=\frac{\sigma^2(\Delta(u)|s-s_1|)+\sigma^2(\Delta(u)|t-t_1|)}{\sigma^2(\Delta(u))}
\frac{\sigma^2(\Delta(u))}{2\sigma^2(u^{1/\beta}\tau^*)}(m^{-\epsilon}_{k,l}(u))^2,
\EQNY
for $(t,s), (t_1,s_1)\in[0,S]^2, (k,l)\in K_u.$

In view of Lemma \ref{PICKANDS}
\BQNY\label{PI_2}
\frac{\pl{\sup_{(t,s)\in I_{k,l}(u)}\overline{Z_u}(s,t)>m^{-\epsilon}_{k,l}(u)}}
{\Psi(m^{-\epsilon}_{k,l}(u))}\to \mathcal{H}_V([0,S]^2)=
\left(\mathcal{H}_{B_{\alpha_\IF}}[0,S]\right)^2,\ \ u\rw\IF,
\EQNY
where $V(t,s)=B^{(1)}_{{\alpha_\IF}}(t)+B^{(2)}_{{\alpha_\IF}}(s)$
with $B_{{\alpha_\IF}}^{(1)}$ and $B^{(2)}_{{\alpha_\IF}}$ being
independent fBms with index ${\alpha_\IF}$.
The rest of the proof follows the same idea as the proof of Theorem \ref{main3}.
\QED

\subsection{Proof of Theorem \ref{main6}}
Similarly to (\ref{F0}), we have
\BQN\label{Last}
\pi^{\inf}_{T_u}(u)\leq \psi^{\inf}_{T_u}(u)\leq \pi^{\inf}_{T_u}(u)+\pl{\inf_{t\in[0,T_uu^{-1/\beta}]}\sup_{s-t\notin E(u)}Z_u(s,t)>m(u)},
\EQN
where
\BQNY
\pi^{\inf}_{T_u}(u)=\pl{\inf_{t\in[0,T_uu^{-1/\beta}]}\sup_{s-t\in E(u)}Z_u(s,t)>m(u)}.
\EQNY
Due to Lemma \ref{L4}, we get
\BQN\label{Last1}
\pl{\inf_{t\in[0,T_uu^{-1/\beta}]}\sup_{s-t\notin E(u)}Z_u(s,t)>m(u)}&\leq& \pl{\sup_{t\in[0,T_uu^{-1/\beta}]}\sup_{s-t\notin E(u)}Z_u(s,t)>m(u)}\nonumber\\
&=&o\left(\frac{u^{1/\beta}}{\Delta(u)m(u)}\Psi(m(u))\right).
\EQN
Next we focus on the asymptotics of $\pi^{\inf}_{T_u}(u)$.
\\
{\it \underline{Case $\varphi=0$ and $\rho\in(0,\IF)$}.}
In order to get the asymptotics of $\pi^{\inf}_{T_u}(u)$
we slightly modify arguments used in
(\ref{U}) and (\ref{L}).
Let
$D(\rho\pm\epsilon,u)=[0, (\rho\pm\epsilon)\frac{\Delta(u)}{u^{1/\beta}}]$ and $F_l(u)=[\tau_u+l\frac{\Delta(u)}{u^{1/\beta}}S,\tau_u+(l+1)\frac{\Delta(u)}{u^{1/\beta}}S]$.
Note that functional $\Phi:=\inf\sup $ satisfies {\bf F1-\bf F2}.
Using that {\bf P1-P4} have been checked in the proof of Theorem 3.1,
following Lemma \ref{PICKANDS}, we have
 \BQNY
\pi^{\inf}_{T_u}(u)
&\leq&
\sum_{l=-N_{S,u}-1}^{N_{S,u}+2}\pl{\inf_{t\in D(\rho+\epsilon,u)}\sup_{s\in F_l(u)}Z_u(s,t)>m(u)}\\
&\leq& \sum_{l=-N_{S,u}-1}^{N_{S,u}+2}\pl{\inf_{t\in D(\rho+\epsilon,u)}\sup_{s\in F_l(u)}\overline{Z_u}(s,t)>m^{-\epsilon}_{0,l}(u)}\\
&\leq& \sum_{l=-N_{S,u}-1}^{N_{S,u}+2}\mathcal{H}_{B_{\alpha_0}}^{\inf}[0,\rho+\epsilon]\mathcal{H}_{B_{\alpha_0}}[0,S]\Psi(m_{0,l}^{-\epsilon}(u))(1+o(1))\\
&\sim& \mathcal{H}_{B_{{\alpha_0}}}^{\inf}[0,\rho+\epsilon]\mathcal{H}_{B_{{\alpha_0}}}(b-\epsilon)^{-1/2}\sqrt{\pi}\frac{u^{1/\beta}}{\Delta(u)m(u)}\Psi(m(u)),
\EQNY
with $N_{S,u}$ and $m^{-\epsilon}_{k,l}(u)$ defined as in the proof of Theorem \ref{main3}.
Similarly,
\BQNY
 \pi^{\inf}_{T_u}(u)\geq\sum_{l=-N_{S,u}}^{N_{S,u}}\pl{\inf_{t\in D(\rho-\epsilon,u)}\sup_{s\in F_l(u)}Z_u(s,t)>m(u)}
 -\sum_{i=1}^2\Sigma_i''(u),
 \EQNY
 with
 \BQNY
  \Sigma_i''(u)&=&\sum_{(l,l_1)\in\Gamma_{i}'}\pl{\inf_{t\in D(\rho-\epsilon,u)}\sup_{s\in F_l(u)}Z_u(s,t)>m(u)\inf_{t\in D(\rho-\epsilon,u)}\sup_{s\in F_{l_1}(u)}Z_u(s_1,t_1)>m(u)}, \ \ i=1,2,
 \EQNY
 where $\Gamma_{1}'=\{(l,l_1), -N_{S,u}\leq l< l_1+1\leq N_{S,u} \}$ and $\Gamma_{2}'=\{(l,l_1), -N_{S,u}\leq l= l_1+1\leq N_{S,u} \}$.\\
Clearly (by the proof of Theorem \ref{main3})
\begin{eqnarray*}
\sum_{i=1}^2\Sigma_i''(u)
&\leq&
\sum_{i=1}^2\sum_{(l,l_1)\in\Gamma_{i}'}
\pl{\sup_{t\in D(\rho-\epsilon,u)}\sup_{s\in F_l(u)}Z_u(s,t)>m(u)\sup_{t\in D(\rho-\epsilon,u)}\sup_{s\in F_{l_1}(u)}Z_u(s_1,t_1)>m(u)}\\
&=&
o\left(\frac{u^{1/\beta}}{\Delta(u)m(u)}\Psi(m(u))\right)
\end{eqnarray*}
and
\[
\sum_{l=-N_{S,u}}^{N_{S,u}}\pl{\inf_{t\in D(\rho-\epsilon,u)}\sup_{s\in F_l(u)}Z_u(s,t)>m(u)}
\geq
\mathcal{H}_{B_{{\alpha_0}}}^{\inf}[0,\rho-\epsilon]\mathcal{H}_{B_{{\alpha_0}}}(b+\epsilon)^{-1/2}\sqrt{\pi}\frac{u^{1/\beta}}{\Delta(u)m(u)}\Psi(m(u)).
\]

Thus, letting $\epsilon\rw 0$, in view of  (\ref{Last}) and (\ref{Last1}), we obtain
\BQN\label{from.3}
\psi^{\inf}_{T_u}(u)\sim \mathcal{H}_{B_{{\alpha_0}}}^{\inf}[0,\rho]\mathcal{H}_{B_{{\alpha_0}}}b^{-1/2}\sqrt{\pi}\frac{u^{1/\beta}}{\Delta(u)m(u)}\Psi(m(u)).
\EQN
\\
{\it \underline{Case $\varphi=0$ and $\rho=0$}}.
The idea of proof is based on the observation that
\BQNY
\psi^{\inf}_{\epsilon\Delta(u)}(u)\leq\psi^{\inf}_{T_u}(u)\leq \psi_0(u)
\EQNY
holds for any $\epsilon>0$ and $u$ sufficiently large.
Following (\ref{from.3}),
$\psi^{\inf}_{\epsilon\Delta(u)}(u)
=\mathcal{H}_{B_{{\alpha_0}}}^{\inf}[0,\epsilon]\mathcal{H}_{B_{{\alpha_0}}}b^{-1/2}\sqrt{\pi}\frac{u^{1/\beta}}{\Delta(u)m(u)}\Psi(m(u))
(1+o(1))$ as $u\to\infty$.
Using that, due to \cite{DI2005},
$\psi_0(u)=
\mathcal{H}_{B_{{\alpha_0}}}b^{-1/2}\sqrt{\pi}\frac{u^{1/\beta}}{\Delta(u)m(u)}\Psi(m(u))
(1+o(1))$ as $u\to\infty$
and
$\lim_{\epsilon\to 0}\mathcal{H}_{B_{{\alpha_0}}}^{\inf}[0,\epsilon]=1$,
the proof is completed.
\\
{\it \underline{Case $\varphi\in(0,\IF]$ with $\rho\in[0,\IF)$}}.
The proof of this case can be established in the same way as presented the above.
\QED

\section{Appendix}\label{ap.a}
In the appendix we present the  proofs of Lemma \ref{PICKANDS}-\ref{L3}.\\
\prooflem{PICKANDS}.
Since in large part the proof is the same as the proof of Lemma 2 in \cite{DI2005}, we present
the steps that confirm extension to the class of  continuous functionals $\Phi$ that satisfy {\bf F1-F2}.
By the classical transformation, for any $\vk{k}_u\in K_u$, we have
\BQN\label{IT}
&&\mathbb{P}\left(\Phi(X^{(u,\vk{k}_u)})>g_{\vk{k}_u}(u)\right)\nonumber\\
&=&\frac{1}{\sqrt{2\pi}g_{\vk{k}_u}(u)}e^{-\frac{1}{2}g_{\vk{k}_u}^2(u)}
\int_{\mathbb{R}}e^we^{-\frac{1}{2}\frac{w^2}{g^2_{\vk{k}_u}(u)}}\mathbb{P}\left(\Phi(X^{(u,\vk{k}_u)})>g_{\vk{k}_u}(u)\Bigl\lvert X_0^{(u,\vk{k}_u)}=g_{\vk{k}_u}(u)-\frac{w}{g_{\vk{k}_u}(u)}\right)dw
\EQN
In light of {\bf F2}, we have
\BQNY
&&\mathbb{P}\left(\Phi(X^{(u,\vk{k}_u)})>g_{\vk{k}_u}(u)\Bigl\lvert X_0^{(u,\vk{k}_u)}=g_{\vk{k}_u}(u)-\frac{w}{g_{\vk{k}_u}(u)}\right)\\
&=&\mathbb{P}\left(\Phi\left(g_{\vk{k}_u}(u)\left(X_{\vk{t}}^{(u,\vk{k}_u)}-r_{u,\vk{k}_u}(t)X_{\vk{0}}^{(u,\vk{k}_u)}\right)
-g_{\vk{k}_u}^2(u)(1-r_{u,\vk{k}_u}(\vk{t}))+w(1-r_{u,\vk{k}_u}(\vk{t}))\right)> w\right),
\EQNY
with $r_{u,\vk{k}_u}(\vk{t})=\mathbb{E}\left(X_{\vk{t}}^{(u,\vk{k}_u)}X_{\vk{0}}^{(u,\vk{k}_u)}\right)$.
The reasoning as used in Lemma 2 in \cite{DI2005}, in view of (\ref{sv1}), {\bf F1-F2} and {\bf P1-P4} , implies that
\BQNY
\Phi\left(g_{\vk{k}_u}(u)\left(X_{\vk{t}}^{(u,\vk{k}_u)}-r_{u,\vk{k}_u}(\vk{t})X_{\vk{0}}^{(u,\vk{k}_u)}\right)
-g_{\vk{k}_u}^2(u)(1-r_{u,\vk{k}_u}(\vk{t}))+w(1-r_{u,\vk{k}_u}(\vk{t}))\right)
\EQNY
weakly converges to
$\Phi\left(\sqrt{2}V(\vk{t})-\sigma^2_V(\vk{t})\right).
$
Besides, (\ref{Cor}) and {\bf P3} lead to, for $u$ large enough,
\BQNY
g_{\vk{k}_u}^2(u)\mathbb{E}\left(X_{\vk{t}}^{(u,\vk{k}_u)}-r_{u,\vk{k}_u}(\vk{t})X_{\vk{0}}^{(u,\vk{k}_u)}-
X_{\vk{s}}^{(u,\vk{k}_u)}+r_{u,\vk{k}_u}(\vk{s})X_{\vk{0}}^{(u,\vk{k}_u)}\right)^2\leq \mathbb{Q}\sum_{i=1}^d|s_i-t_i|^{\eta_i}.
\EQNY
Thus by {\bf F1} and Fernique inequality, we derive for $u$ large enough,
\BQNY
&&\mathbb{P}\left(\Phi\left(g_{\vk{k}_u}(u)\left(X_{\vk{t}}^{(u,\vk{k}_u)}-r_{u,\vk{k}_u}(\vk{t})X_{\vk{0}}^{(u,\vk{k}_u)}\right)
-g_{\vk{k}_u}^2(u)(1-r_{u,\vk{k}_u}(\vk{t}))+w(1-r_{u,\vk{k}_u}(\vk{t}))\right)> w\right)\\
&\leq&\mathbb{P}\left(\sup_{\vk{t}\in M}\left(g_{\vk{k}_u}(u)\left(X_{\vk{t}}^{(u,\vk{k}_u)}-r_{u,\vk{k}_u}(\vk{t})X_{\vk{0}}^{(u,\vk{k}_u)}\right)
-g_{\vk{k}_u}^2(u)(1-r_{u,\vk{k}_u}(\vk{t}))+w(1-r_{u,\vk{k}_u}(\vk{t}))\right)> w\right)\\
&\leq&\mathbb{P}\left(\sup_{\vk{t}\in M}g_{\vk{k}_u}(u)\left(X_{\vk{t}}^{(u,\vk{k}_u)}-r_{u,\vk{k}_u}(\vk{t})X_{\vk{0}}^{(u,\vk{k}_u)}\right)
> w-a_1\right)\\
&\leq& a_2e^{-a_3(w-a_1)^2}
\EQNY
with $a_i, i=1, 2, 3$ positive constants.
The above gives the function that (uniformly) dominates the integrant in (\ref{IT}).
Then using the dominated convergence theorem, we can get the claim.
\QED

\prooflem{R1}. Since the upper bound is straightforward,
we focus on the proof that
$\sigma^2(t)\ge C_1 t^2$ in a neighbourhood of 0.
For this we
use a slight modification of the arguments given in \cite{DE2002}.

From {\bf AI},
there exists $T_0>0$ such that for all $T\ge T_0$ we have
$\sigma(T)>0$ and $\dot{\sigma^2}(T)>0$.

Observe that
\BQNY
\sigma(T_0)\sigma(t)\geq \E{X(T_0)X(t)}\ge 2^{-1}\left(\sigma^2(T_0)-\sigma^2(|T_0-t|)\right).
\EQNY
Thus, by Taylor expansion, with $\rho_t\in (0,t)$ (and $t>0$ small), we get
\BQNY
\sigma^2(T)-\sigma^2(T-t)=\dot{\sigma^2}(T-\rho_t)t\leq 2\sigma(T)\sigma(t),
\EQNY
which implies that
$ \sigma^2(t)\geq \left(\frac{\dot{\sigma^2}(T_0)}{4\sigma(T_0)}\right)^2t^2$
in a neighbourhood of zero.\QED\\

 \prooflem{L1}.
 Recall that
 $\sigma_{u}(\tau)=\frac{\sigma(u^{1/\beta}\tau)}{\sigma(u^{1/\beta})(1+c\tau^\beta)}$.
 By UCT (see, e.g., Theorem 1.5.2 in \cite{BI1989})
 we
 have that
\BQN\label{V1}
\lim_{u\rw\IF}\sigma_{u}(\tau)=\frac{\tau^{\alpha_\IF}}{1+c\tau^\beta}=g(\tau)
\EQN
 holds uniformly on $(0,S]$ for any $S>0$.
 Moreover $\tau^*=\left(\frac{{\alpha_\IF}}{c(\beta-{\alpha_\IF})}\right)^{1/\beta}$
 is the unique maximizer of $g(\tau)$.
Further, by Potter's theorem in (see, e.g., \cite{BI1989}), for any $0<\epsilon<\beta-{\alpha_\IF}$ there exists a constant $u_\epsilon>0$ such that for all $\tau>1$ and $u>u_\epsilon$, we have
\BQN\label{V2}
\sigma_{u}(\tau)\leq \frac{(1+\epsilon)\tau^{{\alpha_\IF}+\epsilon}}{1+c\tau^\beta}\rw 0,
\EQN
as $\tau\rw \IF$.
Combing (\ref{V1}) with (\ref{V2})
we conclude that there exist $S_1, S_2$ such that
for sufficiently large $u$
the maximum of $\sigma_u(\tau)$ is attained in $[S_1, S_2]$ with $0<S_1<\tau^*<S_2<\IF$.
Moreover, by {\bf AI},
\BQN\label{V3}
\dot{\sigma_u}(\tau)\RW \dot{g}(\tau),\ \ \ \ddot{\sigma}_u(\tau)\RW\ddot{g}(\tau),\ \ \ \tau\in[S_1,S_2]
\EQN
and
$\ddot{g}(\tau^*)<0$.

The above implies that, for each sufficiently large $u$,
there exists unique $\tau_u$ such that
$\tau_u\to\tau^*$ as $u\to\infty$,
$\dot{\sigma_u}(\tau_u)=0$
and
$\ddot{\sigma}_u(\tau_u)<0$. This implies that $\tau_u$ is the
unique maximizer of $\sigma_u(\tau)$, for sufficiently large $u$.

It is straightforward to check that
\BQNY
\frac{g(\tau)}{g(\tau^*)}=1-\frac{B}{2A}(\tau-\tau^*)^2(1+o(1)),\ \ \ \tau\rw\tau^*,
\EQNY
which combined with (\ref{V1}) and (\ref{V3}) yields (\ref{V4}). \QED\\
\prooflem{L2}. By direct calculations,
\BQNY
\lefteqn{1-r_u(s,t,s_1,t_1)=}\\
&=&\frac{2\sigma(u^{1/\beta}(s-t))\sigma(u^{1/\beta}(s_1-t_1))+\sigma^2(u^{1/\beta}|t-t_1|)+\sigma^2(u^{1/\beta}|s-s_1|)-
\sigma^2(u^{1/\beta}(s-t_1))-\sigma^2(u^{1/\beta}(s_1-t))}
{2\sigma(u^{1/\beta}(s-t))\sigma(u^{1/\beta}(s_1-t_1))}\\
&=&\frac{D_u^{(1)}(s,t,s_1,t_1)-D_u^{(2)}(s,t,s_1,t_1)+D_u^{(3)}(s,t,s_1,t_1)}
{2\sigma(u^{1/\beta}(s-t))\sigma(u^{1/\beta}(s_1-t_1))},
\EQNY
where
\BQNY
D_u^{(1)}(s,t,s_1,t_1)&=&\sigma^2(u^{1/\beta}|t-t_1|)+\sigma^2(u^{1/\beta}|s-s_1|),\\ D_u^{(2)}(s,t,s_1,t_1)&=&\left(\sigma(u^{1/\beta}(s-t)-\sigma(u^{1/\beta}(s_1-t_1))\right)^2,\\
D_u^{(3)}(s,t,s_1,t_1)&=&\sigma^2(u^{1/\beta}(s-t))+\sigma^2(u^{1/\beta}(s_1-t_1))-\sigma^2(u^{1/\beta}(s-t_1))-\sigma^2(u^{1/\beta}(s_1-t)).
\EQNY
Due to UCT, as $u\rw\IF$,
\BQN\label{V5}
\frac{\sigma^2(u^{1/\beta})t^2}{\sigma^2(u^{1/\beta}t)}
\RW t^{2-2{\alpha_\IF}},\ \ \ t\in(0,S],\ S>0.
\EQN
It follows from mean value theorem and (\ref{V5})
that for $|t-t_1|\leq \delta_u,s-t,s_1-t_1\in E(u)$, with $\theta\in E(u)$,
\BQN\label{V6}
\frac{D_u^{(2)}(s,t,s_1,t_1)}{D_u^{(1)}(s,t,s_1,t_1)}&=&\frac{\left(u^{1/\beta}\dot{\sigma}(u^{1/\beta}\theta)(s-s_1-t+t_1)\right)^2}
{\sigma^2(u^{1/\beta}|t-t_1|)+\sigma^2(u^{1/\beta}|s-s_1|)}\sim
\frac{\alpha_\IF^2\sigma^2(u^{1/\beta}\theta)(s-s_1-t+t_1)^2}{\theta^2\left(\sigma^2(u^{1/\beta}|t-t_1|)+\sigma^2(u^{1/\beta}|s-s_1|)\right)}\\
&\leq&\frac{2\alpha_\IF^2\sigma^2(u^{1/\beta}\theta)((s-s_1)^2+(t-t_1)^2)}{\theta^2\left(\sigma^2(u^{1/\beta}|t-t_1|)+\sigma^2(u^{1/\beta}|s-s_1|)\right)}\nonumber\\
&\leq&\frac{2\alpha_\IF^2\sigma^2(u^{1/\beta}\theta)(s-s_1)^2}{\theta^2\sigma^2(u^{1/\beta}|s-s_1|)}+
\frac{2\alpha_\IF^2\sigma^2(u^{1/\beta}\theta)(t-t_1)^2}{\theta^2\sigma^2(u^{1/\beta}|t-t_1|)}
\rw 0,\ \ \ u\rw \IF.\nonumber
\EQN
Using Taylor expansion, we have
\BQNY
D_u^{(3)}(s,t,s_1,t_1)
&=&u^{1/\beta}\dot{\sigma^2}(u^{1/\beta}(s-t_1))(t_1-t)+\frac{1}{2}u^{2/\beta}\ddot{\sigma^2}(u^{1/\beta}\theta_1)(t-t_1)^2\\
&&+u^{1/\beta}\dot{\sigma^2}(u^{1/\beta}(s_1-t))(t-t_1)
+\frac{1}{2}u^{2/\beta}\ddot{\sigma^2}(u^{1/\beta}\theta_2)(t-t_1)^2\\
&=&\frac{1}{2}u^{2/\beta}\ddot{\sigma^2}(u^{1/\beta}\theta_1)(t-t_1)^2+\frac{1}{2}u^{2/\beta}\ddot{\sigma^2}(u^{1/\beta}\theta_2)(t-t_1)^2\\
&&+u^{2/\beta}\ddot{\sigma^2}(u^{1/\beta}\theta_3)(t_1-t)(s-s_1+t-t_1)\\
&\leq&u^{2/\beta}\left(\frac{1}{2}\ddot{\sigma^2}(u^{1/\beta}\theta_1)+\frac{1}{2}\ddot{\sigma^2}(u^{1/\beta}\theta_2)+2\ddot{\sigma^2}(u^{1/\beta}\theta_3)\right)(t-t_1)^2\\
&&+2u^{2/\beta}\ddot{\sigma^2}(u^{1/\beta}\theta_3)(s-s_1)^2,
\EQNY
where $\theta_1$, $\theta_2$ and $\theta_3$ are some positive constants satisfying $\frac{\tau^*}{2}<\theta_i<\frac{3}{2}\tau^*, i=1,2,3$,  for $u$ sufficiently large.
Similarly, in the light of (\ref{V5}), for $|t-t_1|\leq \delta_u,s-t,s_1-t_1\in E(u)$,
\BQN\label{V7}
\frac{D_u^{(3)}(s,t,s_1,t_1)}{D_u^{(1)}(s,t,s_1,t_1)}\rw 0,\ \ \ u\rw \IF.
\EQN
Hence, the combination of (\ref{V6}) and (\ref{V7})
implies the assertion.\QED\\
 \prooflem{L3}. Substituting $s$ and $s_1$ by $t+\tau$ and $t_1+\tau_1$ respectively yields
\BQNY
\lefteqn{r_u(t+\tau,t,t_1+\tau_1,t_1)}\\
&=&\frac{\sigma^2(u^{1/\beta}|t-t_1+\tau|)+\sigma^2(u^{1/\beta}|t_1-t+\tau_1|)-\sigma^2(u^{1/\beta}|t-t_1+\tau-\tau_1|)
-\sigma^2(u^{1/\beta}|t-t_1|)}{2\sigma(u^{1/\beta}\tau)\sigma(u^{1/\beta}\tau_1)}.
\EQNY
Now suppose $t_1>t$ and $t_1-t>R$ with $R$ a large enough positive constant. Using Taylor expansion at point $t_1-t$, we have
\BQNY
&&\sigma^2(u^{1/\beta}(t_1-t-\tau))+\sigma^2(u^{1/\beta}(t_1-t+\tau_1))-\sigma^2(u^{1/\beta}(t_1-t+\tau_1-\tau))
-\sigma^2(u^{1/\beta}(t_1-t))\\
&=&\sigma^2(u^{1/\beta}(t_1-t))-\dot{\sigma^2}(u^{1/\beta}(t_1-t))u^{1/\beta}\tau+\frac{1}{2}\ddot{\sigma^2}(u^{1/\beta}(t_1-t+\theta_1(u)))
u^{2/\beta}\tau^2\\
&&+\sigma^2(u^{1/\beta}(t_1-t))+\dot{\sigma^2}(u^{1/\beta}(t_1-t))u^{1/\beta}\tau_1+\frac{1}{2}\ddot{\sigma^2}(u^{1/\beta}
(t_1-t+\theta_2(u)))u^{2/\beta}\tau_1^2\\
&&-\left(\sigma^2(u^{1/\beta}(t_1-t))+\dot{\sigma^2}(u^{1/\beta}(t_1-t))u^{1/\beta}(\tau_1-\tau)+\frac{1}{2}\ddot{\sigma^2}(u^{1/\beta}(t_1-t+\theta_3(u)))
u^{2/\beta}(\tau_1-\tau)^2\right)\\
&&-\sigma^2(u^{1/\beta}(t_1-t))\\
&=&\frac{1}{2}\ddot{\sigma^2}(u^{1/\beta}(t_1-t+\theta_1(u)))u^{2/\beta}\tau^2+\frac{1}{2}\ddot{\sigma^2}(u^{1/\beta}
(t_1-t+\theta_2(u)))u^{2/\beta}\tau_1^2\\
&&-\frac{1}{2}\ddot{\sigma^2}(u^{1/\beta}(t_1-t+\theta_3(u))u^{2/\beta}(\tau_1-\tau)^2,
\EQNY
where $\theta_i(u), i=1,2,3$ are some constant satisfying $|\theta_i(u)|\leq 2\tau^*, i=1,2,3$ for $u$ large enough.
Further, by {\bf AI}, we have
\BQNY
\left|\frac{\ddot{\sigma^2}(u^{1/\beta}(t_1-t+\theta))u^{2/\beta}}{\sigma(u^{1/\beta}\tau)\sigma(u^{1/\beta}\tau_1)}\right|
&=&\left|\frac{\ddot{\sigma^2}(u^{1/\beta}(t_1-t+\theta))u^{2/\beta}(t_1-t+\theta)^2}{\sigma^2(u^{1/\beta}(t_1-t+\theta))}
\frac{\sigma^2(u^{1/\beta}(t_1-t+\theta))}{\sigma(u^{1/\beta}\tau)\sigma(u^{1/\beta}\tau_1)(t_1-t+\theta)^2}\right|\\
&\leq & \mathbb{Q}2\alpha_\IF|2\alpha_\IF-1| \frac{(t_1-t+\theta)^{2\alpha_\IF+\epsilon}}{(\tau^*)^{2\alpha_\IF+\epsilon}(t_1-t+\theta)^2}\\
&\leq&\frac{\mathbb{Q}_1}{(t_1-t+\theta)^{2-2\alpha_\IF-\epsilon}},
\EQNY
where $\mathbb{Q}$ and $\mathbb{Q}_1$ are two fixed positive constants, $|\theta|\leq 2\tau^*$ and $0<\epsilon<2-2\alpha_\IF$.
Thus we have, as $R\rw \IF$,
\BQNY
\frac{\ddot{\sigma^2}(u^{1/\beta}(t_1-t+\theta))u^{2/\beta}}{\sigma(u^{1/\beta}\tau)\sigma(u^{1/\beta}\tau_1)}\RW 0,
\EQNY which implies that for $u$ large enough, $|t-t_1|>R ,\tau,\tau_1\in E(u)$
$$r_u(t+\tau,t,t_1+\tau_1,t_1)\RW 0,\ \ \ R\rw \IF.$$
Next we concentrate on the case of $|t-t_1|\leq R ,\tau,\tau_1\in E(u)$ with $R$ a positive constant. Applying UCT, we have
\BQNY
r_u(t+\tau,t,t_1+\tau_1,t_1)&\RW& \frac{|t-t_1+\tau|^{2\alpha_\IF}+|t_1-t+\tau_1|^{2\alpha_\IF}-|t-t_1+\tau-\tau_1|^{2\alpha_\IF}-|t-t_1|^{2\alpha_\IF}}{2\tau^{\alpha_\IF}\tau_1^{\alpha_\IF}}\\
&\RW& \frac{1}{2}\left(|1+\frac{t-t_1}{\sqrt{\tau\tau_1}}|^{2\alpha_\IF}+|1+\frac{t_1-t}{\sqrt{\tau\tau_1}}|^{2\alpha_\IF}-2|\frac{t-t_1}
{\sqrt{\tau\tau_1}}|^{2\alpha_\IF}\right)=f(x),
\EQNY
with $\frac{t-t_1}{\sqrt{\tau\tau_1}}=x$.
It is straightforward to check that $\sup_{x\in[\delta,\IF)}|f(x)|<1$ for any $\delta>0$.
This completes the proof.\QED\\

{\bf Acknowledgements} We are thankful to Enkelejd Hashorva for proposing the topic of this paper and various related discussions. We also kindly acknowledge partial
support from the Swiss National Science Foundation Project 200021-140633/1,
and the project RARE -318984  (an FP7  Marie Curie IRSES Fellowship).
KD also acknowledges partial support by NCN Grant No 2013/09/B/ST1/01778 (2014-2016).

\bibliographystyle{plain}

 \bibliography{storageP}

\def\lfhook#1{\setbox0=\hbox{#1}{\ooalign{\hidewidth
  \lower1.5ex\hbox{'}\hidewidth\crcr\unhbox0}}}
  \def\polhk#1{\setbox0=\hbox{#1}{\ooalign{\hidewidth
  \lower1.5ex\hbox{`}\hidewidth\crcr\unhbox0}}}
  \def\polhk#1{\setbox0=\hbox{#1}{\ooalign{\hidewidth
  \lower1.5ex\hbox{`}\hidewidth\crcr\unhbox0}}} \def\cprime{$'$}
  \def\cprime{$'$}
\begin{thebibliography}{10}

\bibitem{AdlerTaylor}
R.J. Adler and J.E. Taylor.
\newblock {\em Random fields and geometry}.
\newblock Springer Monographs in Mathematics. Springer, New York, 2007.

\bibitem{AS2004}
J.~M.~P. Albin and G.~Samorodnitsky.
\newblock On overload in a storage model, with a self-similar and infinitely
  divisible input.
\newblock {\em Ann. Appl. Probab.}, 14(2):820--844, 2004.

\bibitem{BI1989}
N.~H. Bingham, C.~M. Goldie, and J.~L. Teugels.
\newblock {\em Regular variation}, volume~27 of {\em Encyclopedia of
  Mathematics and its Applications}.
\newblock Cambridge University Press, Cambridge, 1989.

\bibitem{DE2014}
K.~D\c{e}bicki and K.~M. Kosi{\'n}ski.
\newblock On the infimum attained by the reflected fractional {B}rownian
  motion.
\newblock {\em Extremes}, 17(3):431--446, 2014.

\bibitem{DE2002}
K.~D{\c{e}}bicki.
\newblock Ruin probability for {G}aussian integrated processes.
\newblock {\em Stochastic Process. Appl.}, 98(1):151--174, 2002.

\bibitem{DI2005}
A.~B. Dieker.
\newblock Extremes of {G}aussian processes over an infinite horizon.
\newblock {\em Stochastic Process. Appl.}, 115(2):207--248, 2005.

\bibitem{HA2013}
E.~Hashorva, L.~Ji, and V.I. Piterbarg.
\newblock On the supremum of {$\gamma$}-reflected processes with fractional
  {B}rownian motion as input.
\newblock {\em Stochastic Process. Appl.}, 123(11):4111--4127, 2013.

\bibitem{HP99}
J.~H{\"u}sler and V.I. Piterbarg.
\newblock Extremes of a certain class of {G}aussian processes.
\newblock {\em Stochastic Process. Appl.}, 83(2):257--271, 1999.

\bibitem{HP2004}
J.~H{\"u}sler and V.I. Piterbarg.
\newblock Limit theorem for maximum of the storage process with fractional
  {B}rownian motion as input.
\newblock {\em Stochastic Process. Appl.}, 114(2):231--250, 2004.

\bibitem{LE1983}
M.R. Leadbetter, G.~Lindgren, and H.~Rootz{\'e}n.
\newblock {\em Extremes and related properties of random sequences and
  processes}, volume~11.
\newblock Springer Verlag, 1983.

\bibitem{PengLith}
P.~Liu, E.~Hashorva, and L.~Ji.
\newblock On the {$\gamma$}-reflected processes with f{B}m input.
\newblock {\em Lithuanian Math J., in press}, 2015.

\bibitem{Norros94}
I.~Norros.
\newblock A storage model with self-similar input.
\newblock {\em Queueing Systems Theory Appl.}, 16(3-4):387--396, 1994.

\bibitem{Pit96}
V.I. Piterbarg.
\newblock {\em Asymptotic methods in the theory of {G}aussian processes and
  fields}, volume 148 of {\em Translations of Mathematical Monographs}.
\newblock American Mathematical Society, Providence, RI, 1996.

\bibitem{PI}
V.I. Piterbarg.
\newblock Large deviations of a storage process with fractional {B}rownian
  motion as input.
\newblock {\em Extremes}, 4(2):147--164, 2001.

\bibitem{REICH}
Edgar Reich.
\newblock On the integrodifferential equation of {T}ak\'acs. {I}.
\newblock {\em Ann. Math. Statist}, 29:563--570, 1958.

\end{thebibliography}

\end{document}